\documentclass[10pt]{article}
\usepackage{amsmath}
\usepackage{latexsym}
\usepackage{amssymb}
\usepackage{amsfonts}
\usepackage{amsthm}
\title{DIFFERENTIAL CORRESPONDENCES \\  AND CONTROL THEORY}
\author{J.-F. Pommaret \\ CERMICS, Ecole des Ponts ParisTech, France \\
 jean-francois.pommaret@wanadoo.fr \\
ORCID:0000-0003-0907-2601}
\date{  }
\begin{document}
\maketitle

\noindent
{\bf ABSTRACT}    \\

When a differential field $K$ having $n$ commuting derivations is given together with two finitely generated differential extensions $L$ and $M$ of $K$, an important problem in differential algebra is to exhibit a common differential extension $N$ in order to define the new differential extensions $L\cap M$ and the smallest differential field $(L,M)\subset N$ containing both $L$ and $M$. Such a result allows to generalize the use of complex numbers in classical algebra. Having now two finitely generated differential modules $L$ and $M$ over the non-commutative ring ring $D=K[d_1,... ,d_n]=K[d]$ of differential operators with coefficients in $K$, we may similarly look for a differential module $N$ containing both $L$ and $M$ in order to define $L\cap M$ and $L+M$. This is {\it exactly} the situation met in linear or non-linear OD or PD control theory by selecting the inputs and the outputs among the control variables. However, in many recent books and papers, we have shown that controllability was a {\it built-in} property of a control system, not depending on the choice of inputs and outputs. The main purpose of this paper is to revisit control theory by showing the specific importance of the two previous problems and the part plaid by $N$ in both cases for the parametrization of the control system. The essential tool will be the study of {\it differential correspondences}, a modern name for what was called {\it B\"{a}cklund problem} during the last century, namely the study of elimination theory for groups of variables among systems of linear or nonlinear OD or PD equations. The difficulty is to revisit {\it differential homological algebra} by using non-commutative localization. Finally, when $M$ is a $D$-module, this paper is using for the first time the fact that the system $R=hom_K(M,K)$ is a $D$-module for the Spencer operator acting on sections, avoiding thus behaviors, trajectories and signal spaces.  \\

\vspace{3cm}

\noindent
{\bf KEY WORDS}    \\
Differential modules; Differential extensions; Differential elimination; Controllability; \\
K\"{a}hler differentials; B\"{a}cklund problem. \\

 \newpage

 \noindent
 {\bf 1) INTRODUCTION}  \\

The story started in 1970 at Princeton university when the author of tis paper was a visiting student of D.C. Spencer and his colleague J. Wheeler from the nearby physics department set up a 1000 \$ challenge for proving that Einstein equations could be parametrized by potential-like functions like Maxwell equations. It is only in 1995 that he found the negative solution of this chalenge, only paid back one dollar (!) by Wheeler because the relativistic community was (and still is !) convinced about the existence of such a parametrization. Accordingly, such a result can only be found today in books of control theory ([22],[25],[44]). Presenting this result at the algebra seminar of M.-P. Malliavin in the Institut Henri Poincar\'{e}of Paris the same year, he found by chance on display in the library the translation from Japanese of the 1970 Master thesis of M. Kashiwara and discovered the usefulness of differential homological algebra that Spencer {\it never} told him about during his stay in Princeton ([10]). \\  

Meanwhile, U. Oberst from Innsbruck university succeeded applying these new tools to control theory, in particular for studying controllability for multi-dimensional systems with constant coefficients ([14],[15]). However, the reader may discover on the net (www.ricam.oeaw.ac.at/oberst) how difficult it is to communicate with people familiar with analysis but not with the formal methods (jet theory, differential sequences, diagram chasing) when studying systems of ordinary differential (OD) or partial differential (PD) equations. In the meantime, the author had become aware of the new methods (tensor products of rings and fields with derivations) used by A. Bialynicki-Birula in order to study "Differential Galois Theory" ([2],[18]) that are largely superseding the approach of E. Kolchin in classical differential algebra ([9],[11]). \\

A possibility to escape from such a situation was to publish as fast as possible a book presenting for the first time in a self-contained way the non-commutative aspect of double duality for the study of systems having coefficients in a differential field $K$ ([21], Zbl 1079.93001). Of course, the difficulty was to use commutative algebra for the {\it graded} modules in order to study the corresponding {\it filtred} modules, a hard task indeed. A more specific application of theses new tools to mathematical physics (general relativity, gauge theory) allowed the author to justify the many doubts he already had since a long time about the origin and existence of gravitational waves and black holes, but this is out of the scope of this paper (Comparing [1] to [32] and [35] needs no comment).  \\
 
In the second section we shall study the linear framework and in the third section we shall study the nonlinear framework, separating in each situation the differential geometric approach from the differential algebraic approach and providing various motivating examples. Many of the results are given without proofs that can be found in the many books ([17-22],[24],[29],[30]) and recent papers ([24],[27],[28],[32],[33]) that we have published. It does not seem that that the link existing between the Spencer operator and non-commutative localization of Ore domains is known. In any case, this result has never been used for applications to control theory or even mathematical physics.  \\   \\

\noindent
{\bf 2) LINEAR CORRESPONDENCES}  \\
 
\noindent
{\bf 2.1) LINEAR SYSTEMS} \\

If $X$ is a manifold of dimension $n$ with local coordinates $(x)=(x^1, ... ,x^n)$, we denote as usual by $T=T(X)$ the {\it tangent bundle} of $X$, by $T^*=T^*(X)$ the {\it cotangent bundle}, by ${\wedge}^rT^*$ the {\it bundle of r-forms} and by $S_qT^*$ the {\it bundle of q-symmetric tensors}. More generally, let $E, F, \dots$ be {\it vector bundles} over $X$ with local coordinates $(x^i,y^k), (x^i,z^l)),\dots$ for $i=1,...,n$, $k=1,...,m$, $l=1,...,p$ simply denoted by $(x,y)$, $(x,z)$, {\it projection} $\pi:E\rightarrow X:(x,y)\rightarrow (x)$ and changes of coordinates $\bar{x}=\varphi (x), \bar{y}=A(x)y$. We shall denote by $E^*$ the vector bundle obtained by inverting the matrix $A$ of the changes of coordinates , exactly like $T^*$ is obtained from $T$. We denote by $\xi:X\rightarrow E: (x)\rightarrow (x,y=\xi(x))$ a (local) section of $E$, Under a change of coordinates, a section transforms like $\bar{\xi}(\varphi(x))=A(x)\xi(x)$ and the changes of the derivatives can also be obtained with more work. We shall denote by $J_q(E)$ the {\it q-jet bundle} of $E$ with local coordinates $(x^i, y^k, y^k_i, y^k_{ij},...)=(x,y_q)$ called {\it jet coordinates} and sections ${\xi}_q:(x)\rightarrow (x,{\xi}^k(x), {\xi}^k_i(x), {\xi}^k_{ij}(x), ...)=(x,{\xi}_q(x))$ transforming like the sections $j_q(\xi):(x) \rightarrow (x,{\xi}^k(x), {\partial}_i{\xi}^k(x), {\partial}_{ij}{\xi}^k(x), ...)=(x,j_q(\xi)(x))$ where both ${\xi}_q$ and $j_q(\xi)$ are over the section $\xi$ of $E$. For any $q\geq 0$, $J_q(E)$ is a vector bundle over $X$ with projection ${\pi}_q$ while $J_{q+r}(E)$ is a vector bundle over $J_q(E)$ with projection ${\pi}^{q+r}_q, \forall r\geq 0$.\\

\noindent
{\bf DEFINITION  1.A.1}: A {\it linear system} of order $q$ on $E$ is a vector sub-bundle $R_q\subset J_q(E)$ and a {\it solution} of $R_q$ is a section $\xi$ of $E$ such that $j_q(\xi)$ is a section of $R_q$. \\
  
Let $\mu=({\mu}_1,...,{\mu}_n)$ be a multi-index with {\it length} ${\mid}\mu{\mid}={\mu}_1+...+{\mu}_n$, {\it class} $i$ if ${\mu}_1=...={\mu}_{i-1}=0,{\mu}_i\neq 0$ and $\mu +1_i=({\mu}_1,...,{\mu}_{i-1},{\mu}_i +1, {\mu}_{i+1},...,{\mu}_n)$. We set $y_q=\{y^k_{\mu}{\mid} 1\leq k\leq m, 0\leq {\mid}\mu{\mid}\leq q\}$ with $y^k_{\mu}=y^k$ when ${\mid}\mu{\mid}=0$. There is a natural way to distinguish the section ${\xi}_q$ from the section $j_q(\xi)$ by introducing the {\it Spencer operator} $d:J_{q+1}(E)\rightarrow T^*\otimes J_q(E)$ with components $(d{\xi}_{q+1})^k_{\mu,i}(x)={\partial}_i{\xi}^k_{\mu}(x)-{\xi}^k_{\mu+1_i}(x)$. The kernel of $d$ consists of sections such that ${\xi}_{q+1}=j_1({\xi}_q)=...=j_{q+1}(f)$. Finally, if $R_q\subset J_q(E)$ is a {\it system} of order $q$ on $E$ locally defined by linear equations ${\Phi}^{\tau}(x,y_q)\equiv a^{\tau\mu}_k(x)y^k_{\mu}=0$, the $r$-{\it prolongation} $R_{q+r}={\rho}_r(R_q)=J_r(R_q)\cap J_{q+r}(E)\subset J_r(J_q(E))$ is locally defined when $r=1$ by the linear equations ${\Phi}^{\tau}(x,y_q)=0, d_i{\Phi}^{\tau}(x,y_{q+1})\equiv a^{\tau\mu}_k(x)y^k_{\mu+1_i}+{\partial}_ia^{\tau\mu}_k(x)y^k_{\mu}=0$ and has {\it symbol} $g_{q+r}=R_{q+r}\cap S_{q+r}T^*\otimes E\subset J_{q+r}(E)$ if one looks at the top order terms. If ${\xi}_{q+1}\in R_{q+1}$ is over ${\xi}_q\in R_q$, differentiating the identity $a^{\tau\mu}_k(x){\xi}^k_{\mu}(x)\equiv 0$ with respect to $x^i$ and substracting the identity $a^{\tau\mu}_k(x){\xi}^k_{\mu+1_i}(x)+{\partial}_ia^{\tau\mu}_k(x){\xi}^k_{\mu}(x)\equiv 0$, we obtain the identity $a^{\tau\mu}_k(x)({\partial}_i{\xi}^k_{\mu}(x)-{\xi}^k_{\mu+1_i}(x))\equiv 0$ and thus the restriction $d:R_{q+1}\rightarrow T^*\otimes R_q$. More generally, we have the restriction:   \\
\[   d: {\wedge}^sT^* \otimes R_{q+1} \rightarrow {\wedge}^{s+1}T^* \otimes R_q: ({\xi}^k_{\mu,I}(x)dx^I) \rightarrow (({\partial}_if^k_{\mu,I}(x) - {\xi}^k_{\mu + 1_i,I}(x))dx^i \wedge dx^I)     \eqno{(1)}\] 
using standard multi-index notation for exterior forms, namely $I=\{i_1 <i_2< ... < i_r\}$, $dx^I=dx^{i_1}\wedge ... \wedge dx^{i_r} \in {\wedge}^rT^*$ for a finite basis, and one can easily check that $d\circ d=0$. The restriction of $-d$ to the symbol is called the {\it Spencer map} $\delta: {\wedge}^sT^*\otimes g_{q+1} \rightarrow {\wedge}^{s+1}T^* \otimes g_q $ and $\delta \circ \delta=0$ similarly, leading to the purely algebraic {\it $\delta$-cohomology} $H^s_{q+r}(g_q)$ ([17]-[20],[22],[30],[43]). \\
    
\noindent
{\bf DEFINITION 1.A.2}: A system $R_q$ is said to be {\it formally integrable} when all the equations of order $q+r$ are obtained by $r$ prolongations {\it only}, $\forall r\geq 0$ or, equivalently, when the projections ${\pi}^{q+r+s}_{q+r}:R_{q+r+s}\rightarrow R^{(s)}_{q+r} \subseteq R_{q+r}$ are epimorphisms $\forall r,s\geq0$.  \\

Finding an intrinsic test has been achieved by D.C. Spencer in 1965 ([43]) along coordinate dependent lines sketched by M. Janet in 1920 ([6]) and providing a {\it Pommaret basis} by using the {\it Janet tabular} in a particular coordinate system called $\delta$-{\it regular} like in  ([17]-[20]).       \\ 

\noindent
{\bf THEOREM 1.A 3}: $R_q$ is {\it formally integrable} ({\it involutive}) if ${\pi}^{q+1}_q:R_{q+1} \rightarrow R_q$ is an epimorphism and $g_q$ is $2$-acyclic with $H^2_{q+r}(g_q)=0, \forall r\geq 0$ (involutive with $H^s_{q+r}(g_q)=0, \forall r\geq 0, \forall s=1,...,n$). When $R_q$ is involutive, there exist $n$ integers ${\alpha}^1_q\geq ... \geq {\alpha}^n_q= \alpha \geq 0$ called {\it characters} and we have $ dim(g_{q+r})= {\Sigma}^{n}_{i=1} \frac{(r+i-1)!}{r! (i-1)!}{\alpha}^i_q$, in particular $dim(g_q)={\alpha}^1_q + ... + {\alpha}^n_q, \,\,\, dim(g_{q+1})={\alpha}^1_q + ... + n {\alpha}^n_q$.  \\

\noindent
{\bf REMARK 1.A.4}: As long as the {\it Prolongation}/{\it Projection} (PP) procedure has not been achieved in order to get an involutive system $R^{(s)}_{q+r}$ for $r,s$ large enough, {\it nothing} can be said about the CC (Fine examples can be found in [30] and the recent [32]). \\

When $R_q$ is involutive, the linear differential operator ${\cal{D}}:E\stackrel{j_q}{\rightarrow} J_q(E)\stackrel{\Phi}{\rightarrow} J_q(E)/R_q=F_0$ of order $q$ is said to be {\it involutive}. Introducing the set of solutions $\Theta \subseteq E$ and the {\it Janet bundles}:  \\
\[ F_r= {\wedge}^rT^*\otimes J_q(E)/({\wedge}^rT^*\otimes R_q + \delta (S_{q+1}T^*\otimes E))   \eqno{(2)}   \] 
we obtain the canonical linear {\it Janet sequence} (Introduced in [17], p 185 + p 391):\\
\[  0 \longrightarrow  \Theta \longrightarrow E \stackrel{\cal{D}}{\longrightarrow} F_0 \stackrel{{\cal{D}}_1}{\longrightarrow}F_1 \stackrel{{\cal{D}}_2}{\longrightarrow} ... \stackrel{{\cal{D}}_n}{\longrightarrow} F_n \longrightarrow 0   \eqno{(3)}     \]
where each other operator, induced by the {\it Spencer operator}, is first order involutive and generates the {\it compatibility conditions} (CC) of the preceding one. Similarly, introducing the {\it Spencer bundles}:  \\ \[   C_r= {\wedge}^rT^* \otimes R_q / \delta ({\wedge}^{r-1}T^* \otimes g_{q+1})   \eqno{(4)}  \] 
we obtain the canonical linear $Spencer$ sequence also induced by the {\it Spencer operator}:  \\
\[   0\longrightarrow \Theta \stackrel{j_q}{\longrightarrow} C_0 \stackrel{D_1}{\longrightarrow} C_1 \stackrel{D_2}{\longrightarrow}... \stackrel{D_n}{\longrightarrow} C_n \longrightarrow  0    \eqno{(5)}   \]

 In the case of analytic systems, the following theorem providing the {\it Cartan-K\"{a}hler} CK) {\it data} is well known though its link with involution is rarely quoted because it is usually presented within the framework of exterior calculus ([17],[26]):  \\

\noindent
{\bf THEOREM 1.A.5} (Cartan-K\"{a}hler): If $R_q\subset J_q(E)$ is a linear involutive and analytic system of order $q$ on $E$, there exists one analytic solution $y^k=f^k(x)$ and only one such that:  \\
1) $(x_0,{\partial}_{\mu}f^k(x_0))$ with $0 \leq \mid \mu \mid\leq q-1$ is a point of $R_{q-1}={\pi}^q_{q-1}(R_q)\subset J_{q-1}(E)$.  \\
2) For $i=1,...,n$ the ${\alpha}^i_q$ parametric derivatives ${\partial}_{\mu}f^k(x)$ of class $i$ are equal for $x^{i+1}=x^{i+1}_0,...,x^n=x^n_0$ to ${\alpha}^i_q$ given analytic functions of $x^1,...,x^i$.  \\
 
 The monomorphism $0\rightarrow J_{q+1}(E) \rightarrow J_1(J_q(E))$ allows to identify $R_{q+1} $ with its image ${\bar{R}}_1$ in $ J_1(R_q)$ and we just need to set $R_q=\bar{E}$ in order to obtain the first order system (Spencer form) ${\bar{R}}_1\subset J_1(\bar{E})$ which is also involutive and analytic while ${\pi}^1_0:{\bar{R}}_1\rightarrow \bar{E}$ is an epimorphism. Studying the respective symbols, we may identify $g_{q+r}$ and ${\bar{g}}_r$ while ${\bar{g}}_1$ is involutive. Looking at the Janet board of multiplicative variables we have ${\bar{\alpha}}^i_1 + {\bar{\beta}}^i_1 = \bar{m}=dim(\bar{E})$ and:  \\
 \[   {\bar{\alpha}}^i_1={\alpha}^i_q+ ... + {\alpha}^n_q={\alpha}^i_{q+1} \Rightarrow {\alpha}^i_q={\bar{\alpha}}^i_1-{\bar{\alpha}}^{i+1}_1={\bar{\beta}}^{i+1}_1-{\bar{\beta}}^i_1  \]
We obtain therefore:  \\

\noindent
{\bf COROLLARY 1.A.6}: If $R_1\subset J_1(E)$ is a first order linear involutive and analytic system such that ${\pi}^1_0:R_1 \rightarrow E$ is an epimorphism, then there exists one analytic solution $y^k=f^k(x)$ and only one, such that:  \\
1) $f^1(x),...,f^{{\beta}^1_1}(x)$ are equal to ${\beta}^1_1$ given constants when $x=x_0$.  \\
2) $f^{{\beta}^i_1+1}(x),..., f^{{\beta}^{i+1}_1}(x)$ are equal to ${\beta}^{i+1}_1-{\beta}^i_1$ given analytic functions of $x^1,...,x^i$ when $x^{i+1}=x^{i+1}_0, ...,x^n=x^n_0$.  \\
3) $f^{{\beta}^n_1+1}(x),...,f^m(x)$ are $m-{\beta}^n_1$ given analytic functions of $x^1,...,x^n$.  \\   \\

\noindent
{\bf 2.2) DIFFERENTIAL MODULES}\\

If $A$ is an associative ring with unit $1 \in A$, a subset $S \subset A$ is called a {\it multiplicative subset} if $1 \in S, 0 \notin S, st\in S, \forall s,t\in S$.. In the commutative case, these conditions are sufficient to localize $A$ at $S$ by constructing the new {\it ring of fractions} $S^{-1}A$ over $A$.  For simplicity, we shall suppose that $A$ is an integral domain (no divisor of zero) and we shall choose $S =A - \{0\}$ in order to introduce the {\it field of fractions} $Q(A)=S^{-1}A=AS^{-1}$. The idea is to exhibit new quantities written $\frac{a}{s}$ with the standard rules:  \\
\[ b (\frac{a}{s})= \frac{ab}{s}, \,\,\, \frac{a}{s} + \frac{b}{t}= \frac{at+bs}{st}, \,\,\,  \frac{a}{s} \frac{b}{t}=\frac{ab}{st} \hspace{2cm}  \forall a,b\in A, \forall s,t \in S\]
The same definition can be used for any module $M$ over $A$ in order to introduce the {\it module of fractions} $S^{-1}M$ over $S^{-1}A$ with the rules:  \\
\[    \frac{a}{s} \frac{x}{t}=\frac{ax}{st}, \,\,\, \frac{x}{s} + \frac{y}{t}= \frac{tx+sy}{st}, \hspace{2cm}  \forall x,y\in M  \]  

\noindent
{\bf DEFINITION 2.2.1}: $t(M)=t_S(M)=\{ x \in M \mid \exists s\in S, sx=0 $ is called the {\it torsion submodule} of $M$ over $S$ and we have the exact sequence 
$0 \rightarrow t(M) \rightarrow M \stackrel{\theta}{\rightarrow} S^{-1}M$ of modules over $A$ where the morphism $\theta$ on the right is $x \rightarrow \frac{x}{1}=\frac{sx}{s}$ and we have $S^{-1}M=S^{-1}A {\otimes}_AM$.  \\

In the non-commutative case considered through all this paper, we shall meet four problems:  \\
$\bullet$  \,\, How to compare $s^{-1}a$ with $as^{-1}$ ?.  \\
$\bullet$ \,\, How to decide when we shall say that $ s^{-1}a=t^{-1}b$ ?.  \\
$\bullet$ \,\, How to multiply $s^{-1}a$ by $t^{-1}b$ ?.  \\
$\bullet$ \,\, How to find a common denominator for $s^{-1}a + t^{-1}b$ ?.   \\

\noindent
{\bf LEMMA 2.2.2}: If there exists a left localization of a noetherian $A$ with respect to $S$, then we must have $Sa \cap As \neq \emptyset$. It follows that $As \cap At\cap S \neq \emptyset$ and two fractions can be multiplied or brought to the same denominator. Finally, $t(M)$ is a submodule of $M$.   \\
                 \\
\noindent
{\it Proof}: Roughly, any right fraction $ a s^{-1}$ can be written as a left fraction $t^{-1}b$, that is we must have $ta=bs$. Now, if we have two fractions $s^{-1} a$ and $t^{-1}b$, we can find $u,v \in A $ such that $us=vt \in S$. Hence, we obtain $s^{-1}a=(us)^{-1} (ua)$ and $t^{-1}b= (vt)^{-1}vb=(us)^{-1} vb$. As for the multiplication 
of fractions, we have $(s^{-1}a)(t^{-1}b)= s^{-1}(at^{-1})b= s^{-1}(u^{-1}c)b=(us)^{-1}(cb)$.  \\
Finally, given $x,y \in t(M)$, we can find $s,t \in S$ such that $sx=0, ty=0$. We may thus find $u,v \in A $ such that $us=vt \in S$ and we get $ us(x+y)=usx + vty=0 \Rightarrow x+y \in t(M)$. Also, we can use $ta=bs$ in order to obtain $t(ax)=(ta)x=(bs)x=b(sx)=0  \Rightarrow ax \in t(M)$.   \\
\hspace*{13cm}       $\Box$    \\

Let $K$ be a {\it differential field} with $n$ commuting derivations $({\partial}_1,...,{\partial}_n)$ and consider the ring $D=K[d_1,...,d_n]=K[d]$ of differential operators with coefficients in $K$ with $n$ commuting formal derivatives satisfying $d_ia=ad_i + {\partial}_ia$ in the operator sense. If $P=a^{\mu}d_{\mu}\in D=K[d]$, the highest value of ${\mid}\mu {\mid}$ with $a^{\mu}\neq 0$ is called the {\it order} of the {\it operator} $P$ and the ring $D$ with multiplication $(P,Q)\longrightarrow P\circ Q=PQ$ is filtred by the order $q$ of the operators. We have the {\it filtration} $0\subset K=D_0\subset D_1\subset  ... \subset D_q \subset ... \subset D_{\infty}=D$. As an algebra, $D$ is generated by $K=D_0$ and $T=D_1/D_0$ with $D_1=K\oplus T$ if we identify an element $\xi={\xi}^id_i\in T$ with the vector field $\xi={\xi}^i(x){\partial}_i$ of differential geometry, but with ${\xi}^i\in K$ now. It follows that $D={ }_DD_D$ is a {\it bimodule} over itself, being at the same time a left $D$-module by the composition $P \longrightarrow QP$ and a right $D$-module by the composition $P \longrightarrow PQ$. We define the {\it adjoint} functor $ad:D \longrightarrow D^{op}:P=a^{\mu}d_{\mu} \longrightarrow  ad(P)=(-1)^{\mid \mu \mid}d_{\mu}a^{\mu}$ and we have $ad(ad(P))=P$ both with $ad(PQ)=ad(Q)ad(P), \forall P,Q\in D$. Such a definition can be extended to any matrix of operators by using the transposed matrix of adjoint operators (See [20]-[22],[23],[24],[31],[33],[36],[37] for more details and applications to control theory or mathematical physics). \\  

\noindent
{\bf PROPOSITION 2.2.3}: $D$ is an Ore domain and $S=D - \{ 0 \} \Rightarrow S^{-1}D = D S^{-1}$.  \\

\noindent
{\it Proof}: Let $U \in S$ and $P \in D$ be given. In order to prove the Ore property for $D$, we must find $V \in S$ and $Q \in D$ such that $VP=QU$. Considering the system $Py=u, Uy=v$, it defines a differential module $M$ over $D$ with the finite presentation $D^2 \rightarrow D \rightarrow M \rightarrow 0$. Now, as we have only one unknown and $D=Dy$ in this sequence, then $M$ is a torsion module and $rk_D(M)=0$. From the additivity property of the differential ranks, if there should be no compatibilty condition (CC), see the example below), then the first morphism on the left should be a monomorphism, a result leading to the contradiction $2 - 1 + 0=0$. Accordingly we can find the operators $ V$ and $Q$ such that $PU^{-1}=V^{-1}Q$. Conversely, if now $V$ and $Q$ are given, using the adjoint functor and the fact that $ad(ad(P))=P, \forall P \in D$, we may obtain $ad(V)$ and $ad(Q)$ such that $ad(P)ad(V)=ad(U)ad(Q)$ as before and thus $V=ad(ad(V))$ and $Q=ad(ad(Q))$ such that $VP=QU$, a result showing that $V^{-1}Q=PU^{-1}$. \\
\hspace*{13cm}   $\Box$  \\

\noindent
{\bf EXAMPLE 2.2.4}: With $m=1, n=2, q=1, K=\mathbb{Q}(x^1,x^2)$ let us consider the two first order operators $U=d_2 \in S, P=d_1 + x^2$. Considering the formal system $d_1y + x^2 y=u, d_2y=v$, we obtain $y=d_2u - d_1v - x^2v$ and thus the involutive system with jet notations:
\[   \left \{ \begin{array}{lcl}
y_2 & = & v  \\
y_1 + x^2 y & =  & u   \\
y  & = & u_2 - v_1 - x^2v
\end{array} \right. \fbox { $ \begin{array}{cc}
 1 & 2 \\
 1 & \bullet  \\
 \bullet & \bullet
 \end{array} $  }  \]
Among the three CC that should exist, only two are non-trivial and provide the new {\it second order} (care !) involutive sytem:  \\
\[   \left \{  \begin{array}{lcl}
A & \equiv & u_{22} - v_{12} - x^2 v_2 - 2v =0  \\
B & \equiv & u_{12} + x^2 u_2 - u - v_{11} - 2 x^2 v_1 - (x^2)^2 v = 0
\end{array} \right. \fbox  { $ \begin{array}{cc}
1 & 2 \\
 1  & \bullet
\end{array} $  }  \]
with the unexpected single first order CC $  C \equiv d_2 B - d_1 A - x^2 A = 0$. We obain therefore the two operator identities:  \\
\[   d_{22} (d_1 + x^2 ) = (d_{12} + x^2 d_2 + 2) d_2,   \,\,\,\,  (d_{12} + x^2 d_2 - 1)(d_1 + x^2)= (d_{11} + 2 x^2 d_1 + (x^2)^2)d_2    \]
leading again to the two unexpected localizations:  \\
\[   (d_1 + x^2) (d_2)^{-1} = (d_{22})^{-1}(d_{12} + x^2 d_2 +2 ) = (d_{12} + x^2 d_2 - 1)^{-1}( d_{11} + 2 x^2 d_1 + (x^2)^2)    \]
Taking the adjoint operators, we get in particular $ (d_1 - x^2) d_{22} = d_2 ( d_{12} - x^2 d_2 + 1)$.  \\
In order to achieve this example and explain why such methods, up to our knowledge, have never been used for applications, it just remains to explain the equality of these two fractions in this framework. Indeed, we obtain easily the unique operator identity $ (d_1 + x_2)d_{22} = d_2 ( d_{12} + x^2 d_2 - 1)$ provided by the last CC  
$d_2B = (d_1 + x^2 )A$ for $u$. Reducing to the same denominator can be done if we use the operator identity:  \\
\[    (d_1 + x^2)(d_{12} + x^2 d_2 + 2) =  d_2(d_{11} + 2 x^2 d_1 + (x^2)^2)   \]
produced by the same last CC $d_2B = (d_1 + x^2 )A$ for $v$. We conclude this example exhibiting the corresponding long exact sequence of differential modules:  \\
\[ 0 \rightarrow D \rightarrow D^2 \rightarrow D^2 \rightarrow D \rightarrow M \rightarrow 0  \]
where we have successively from left to right: $D=DC, D^2 = DA + DB, D^2 = Du + Dv, D=Dy$ with Euler-Poincar\'{e} characteristic $1 - 2 + 2 - 1 =rk_D(M)=0$ because $m=1$.  \\

Accordingly, if $y=(y^1, ... ,y^m)$ are differential indeterminates, then $D$ acts on $y^k$ by setting $d_iy^k=y^k_i \longrightarrow d_{\mu}y^k=y^k_{\mu}$ with $d_iy^k_{\mu}=y^k_{\mu+1_i}$ and $y^k_0=y^k$. We may therefore use the jet coordinates in a formal way as in the previous section. Therefore, if a system of OD/PD equations is written in the form ${\Phi}^{\tau}\equiv a^{\tau\mu}_ky^k_{\mu}=0$ with coefficients $a\in K$, we may introduce the free {\it left} differential module $Dy=Dy^1+ ... +Dy^m\simeq D^m$ and consider the differential {\it module of equations} $I=D\Phi\subset Dy$, both with the residual {\it left differential module} $M=Dy/D\Phi$ or $D$-module and we may set $M={ }_DM \in mod(D)$ if we want to specify the action of the ring of differential operators. We may introduce the formal 
{\it prolongation} with respect to $d_i$ by setting $d_i{\Phi}^{\tau}\equiv a^{\tau\mu}_ky^k_{\mu+1_i}+({\partial}_ia^{\tau\mu}_k)y^k_{\mu}$ in order to induce maps $d_i:M \longrightarrow M:{\bar{y} }^k_{\mu} \longrightarrow {\bar{y}}^k_{\mu+1_i}$ by residue with respect to $I$ if we use to denote the residue $Dy \longrightarrow M: y^k \longrightarrow {\bar{y}}^k$ by a bar like in algebraic geometry. However, for simplicity, we shall not write down the bar when the background will indicate clearly if we are in $Dy$ or in $M$. As a byproduct, the differential modules we shall consider will always be {\it finitely generated} ($k=1,...,m<\infty$) and {\it finitely presented} ($\tau=1, ... ,p<\infty$). Equivalently, introducing the {\it matrix of operators} ${\cal{D}}=(a^{\tau\mu}_kd_{\mu})$ with $m$ columns and $p$ rows, we may introduce the morphism $D^p \stackrel{{\cal{D}}}{\longrightarrow} D^m:(P_{\tau}) \longrightarrow (P_{\tau}{\Phi}^{\tau})$ over $D$ by acting with $D$ {\it on the left of these row vectors} while acting with ${\cal{D}}$ {\it on the right of these row vectors} by composition of operators with $im({\cal{D}})=I$. The {\it presentation} of $M$ is defined by the exact cokernel sequence $D^p \stackrel{{\cal{D}}}{\longrightarrow} D^m \stackrel{p}{\longrightarrow} M \longrightarrow 0 $. We notice that the presentation only depends on $K, D$ and $\Phi$ or $ \cal{D}$, that is to say never refers to the concept of (explicit local or formal) solutions. It follows from its definition that $M$ can be endowed with a {\it quotient filtration} obtained from that of $D^m$ which is defined by the order of the jet coordinates $y_q$ in $D_qy$. We have therefore the {\it inductive limit} $0 \subseteq M_0 \subseteq M_1 \subseteq ... \subseteq M_q \subseteq ... \subseteq M_{\infty}=M$ with $d_iM_q\subseteq M_{q+1}$ and $M=DM_q$ for $q\gg 0$ with prolongations $D_rM_q\subseteq M_{q+r}, \forall q,r\geq 0$. \\

\noindent
{\bf DEFINITION 2.2.5}: An exact sequence of morphisms finishing at $M$ is said to be a {\it resolution} of $M$. If the differential modules involved apart from $M$ are free, that is isomorphic to a certain power of $D$, we shall say that we have a {\it free resolution} of $M$. In the general situation of a sequence $ M' \stackrel{f}{\longrightarrow} M  \stackrel{g}{\longrightarrow} M" $ of modules which may not be exact, we may define the {\it coboundary}, {\it cocycle} and {\it cohomology} at $M$ by setting respectively 
$B =im(f) \subseteq Z = ker(g) \Rightarrow   H=Z/B$ and apply the above result to the various short exact sequences like $ 0 \rightarrow Z \rightarrow M \rightarrow im(g) \rightarrow 0$ or $0 \rightarrow B \rightarrow Z \rightarrow H \rightarrow 0$. The {\it deleted complex}is obtained by replacing $M$ by $0$. Applying $hom_D(\bullet, D)$, we obtaine a sequence that may not be exact. The corresponding cohomology modules, called {\it extension modules} $ ext^i _D(M,D)=ext^i(M)$, are torsion modules for $i\geq 1$, do not depend on the resolution of $M$ and only depend on $K, D$ and $M$ ([3],[10],[21],[36],[38],[40]).    \\

Having in mind that $K$ is a left $D$-module with the action $(D,K) \rightarrow K:(P,a) \rightarrow P(a)$ defined by $(b,a) \rightarrow ba=ab, (d_i,a)\rightarrow d_i(a)={\partial}_ia$ and that $D$ is a bimodule over itself for the composition law of operators, {\it we have only two possible constructions} only depending on $K, D$ and $M$:  \\

\noindent
{\bf DEFINITION 2.2.6}: We may define the {\it right} (care !) differential module $hom_D(M,D)$, using the bimodule structure of ${ }_DD_D$ and setting $(fP)(m)=f(m)P$ while checking that:  \\
\[ ((fP)Q)(m)=((fP)(m))Q=(f(m)P)Q=f(m)PQ=(fPQ)(m) \]

\noindent
{\bf REMARK 2.2.7}: When $M$ admits the finite presentation $D^p \stackrel{{\cal{D}}}{\longrightarrow} D^m \stackrel{p}{\longrightarrow} M \rightarrow 0$, applying $hom_D(\bullet , D)$ we obtain the following long exact sequence of {\it right} (care !) differential modules:   \\
\[  0 \leftarrow  N_D  \longleftarrow D^p_D \longleftarrow  D^m_D \longleftarrow hom_D(M,D) \leftarrow 0  \]
and the so-called {\it Malgrange isomorphism} just amounts to the fact that $ext^0(M)=hom_D(M,D)$. We are immediately facing one of the most delicate problems of this section when dealing  with applications and/or effective computations, a problem not solved in the corresponding literature which has been almost entirely using fields of constants ([14],[15],[41],[42]) though the solution is known since a long time ([3],[20],[38]). Indeed, apart for purely mathematical reasons, the only differential modules to be met are left differential modules. By chance, one has the following theorem describing the functorial {\it  side changing procedure} amounting to replace ${\cal{D}}$ by is formal adjoint $ad({\cal{D}})$. In the differential geometric framework, such a procedure amounts to replace an operator ${\cal{D}}:E \rightarrow F$ by its formal adjoint $ad({\cal{D}}): {\wedge}^nT^* \otimes F^* \rightarrow {\wedge}^nT^* \otimes E^*$ and one may set $ad(E)={\wedge}^nT^* \otimes E^*$ for any vector bundle $E$ over $X$ with $dim(E) = dim(ad(E))$ while reversing the arrows. However, as the formal adjoint of an involutive operator may not even be formally integrable, the formal adjoint of an exact (Janet, Spencer) sequence may not be an exact (Janet, Spencer) sequence at all and this is the motivation for introducing the extension modules. The simplest example when $n=1$ can be found in the study of the double pendulum ([36]). In the time-varying case with $x=t$ and $d=d_t$, the Kalman-type system $dy=A(x)y + B(x)u$ is controllable if and only if the operator $\lambda \rightarrow (d\lambda + \lambda A, \lambda B )$ is injective, that is to say if and only if the matrix $(B, (d-A)B, (d-A)^2B, ... )$ has maximum rank ( Compare to [45] where the adjoint is missing).  \\

\noindent
{\bf THEOREM 2.2.8}: There exists an isomorphism $N_D  \rightarrow N = { }_DN = hom_K({\wedge}^nT^*,N_D) $ with inverse $M={ }_DM \rightarrow M_D= {\wedge}^nT^* {\otimes}_K  M$ .  \\

\noindent
{\it Proof}: First of all, we prove that ${\wedge}^nT^*$ has a natural right module structure over $D$ by introducing the basic volume $n$-form $\alpha = dx^1\wedge \dots \wedge dx^n = dx$ and defining 
$\alpha . P = ad(P)(1)dx, \forall P \in D$. We have $\alpha. \xi = - {\partial}_i {\xi}^idx = - div(\xi) dx = - {\cal{L}}(\xi) \alpha, \forall \xi={\xi}^i  d_i \in T$ where ${\cal{L}}$ is the classical Lie derivative on forms and obtain therefore:  \\
\[ \alpha.(a \xi)= (a \,dx).\xi= - {\partial}_i(a{\xi}^i) dx= - a {\partial}_i{\xi}^i dx- \xi (a)dx= - a {\partial}_i{\xi}^i dx - \alpha.(\xi ( a))   \]
\[ \alpha.(\xi a) =  (\alpha. \xi) a= - a {\partial}_i{\xi}^i  dx  \hspace{1cm} \Rightarrow \hspace{1cm}
\alpha. ({\xi \,a}) = \alpha. (a \xi + \xi(a))   \]
From well known properties of the Lie derivative, we have also:  \\
\[      \alpha .(\xi \eta - \eta \xi) = - ({\cal{L}} (\xi){\cal{L}}(\eta) - {\cal{L}}(\eta){\cal{L}}(\xi)) \alpha = 
- {\cal{L}}([\xi, \eta]) \alpha     \] 
Now, using the adjoint map $ad:D \rightarrow D :  P   \rightarrow ad(P)$, we may introduce the {\it  adjoint functor} $ad: mod(D) \rightarrow mod(D^{op}): M \rightarrow ad(M)$ with $m.P=ad(P) m, \forall m\in M, \forall P \in D$ and we have $m.(PQ)= ad(PQ)m=ad(Q)ad(P)m=(m.P).Q, \forall P,Q \in D$.\\
It remains to introduce the $K$-linear isomorphism $ad(M)  \simeq M_D: m \rightarrow m \otimes \alpha$ with: 
\[  (m\otimes \alpha)a=am \otimes \alpha = m \otimes a \alpha, \hspace{1cm}
(m \otimes \alpha) \xi = - \xi m \otimes \alpha + m \otimes \alpha. \xi ,  \,\,\, \forall m \in M, \forall \xi \in T \] 
and check that $(m\otimes \alpha)(\xi a)= (m\otimes \alpha) (a \xi + \xi(a))$.
These definition are coherent because, when $d$ is any $d_i$, we have $div(d)=0$ and thus 
$(m\otimes \alpha)d = - dm \otimes \alpha$, a result leading to the formula:  \\
\[  (m \otimes \alpha)P=(m \otimes \alpha) (a^{\mu}d_{\mu})=(a^{\mu}m\otimes \alpha)d_{\mu}= 
   (-1)^{\mid \mu \mid}d_{\mu}a^{\mu}m \otimes \alpha = ad(P)m \otimes \alpha  \]
The isomorphism $ad(D) \simeq D_D:P \rightarrow ad(P)$ is also {\it  right} $D$-linear because we have successively: \\
\[   P.Q \rightarrow ad(P.Q)=ad(ad(Q)P)=ad(P)ad(ad(Q))=ad(P)Q  \]
These unexpected results explain why the formal adjoint cannot be avoided in non-commutative localization and is so important for applications ranging from control theory ([21],[22],[27]) to continuum mechanics ([5],[23]) or electromagnetism ([31])  and even general relativity ([33],[37]).   \\
\hspace*{13cm}       $\Box$  \\

\noindent
{\bf DEFINITION 2.2.9}: We define the {\it system} $R=hom_K(M,K)$ and set $R_q=hom_K(M_q,K)$ as the {\it system of order} $q$. We have the {\it projective limit} $R=R_{\infty} \longrightarrow ... \longrightarrow R_q \longrightarrow ... \longrightarrow R_1 \longrightarrow R_0$. It follows that $f_q\in R_q:y^k_{\mu} \longrightarrow f^k_{\mu}\in K$ with $a^{\tau\mu}_kf^k_{\mu}=0$ defines a {\it section at order} $q$ and we may set $f_{\infty}=f\in R$ for a {\it section} of $R$. For an arbitrary differential field $K$, {\it such a definition has nothing to do with the concept of a formal power series solution} ({\it care}).\\

Similarly to the preceding definition, we may define the {\it left} (care !) differential module $hom_K(D,K)$, using again the bimodule structure of $D$ and setting $(Qf)(P)=f(PQ), \forall P,Q \in D$, in particular with 
$  (\xi f)(P)f(P \xi), \forall \xi \in T, \forall P \in D$. However, we should have $(af)(P)=f(Pa) \neq f(aP)=a(f(P)),\forall a \in $, unless $K$ is a field of constants like in most of the literature ([14],[15],[42]).    \\

\noindent
{\bf PROPOSITION 2.2.10}: When $M$ is a left $D$-module, then $R$ is also a left $D$-module. \\

\noindent
{\it Proof}: As $D$ is generated by $K$ and $T$ as we already said, let us define:  \\
\[  (af)(m)=af(m) = f(am), \hspace{5mm} \forall a\in K, \forall m\in M \]
\[ (\xi f)(m)=\xi f(m)-f(\xi m), \hspace{5mm} \forall \xi=a^id_i\in T,\forall m\in M  \]
In the operator sense, it is easy to check that $d_ia=ad_i+{\partial}_ia$ and that $\xi\eta - \eta\xi=[\xi,\eta]$ is the standard bracket of vector fields. Using simply $\partial$ in place of any ${\partial}_i$ and $d$ in place of any $d_i$, we have:  \\
\[  \begin{array}{rcl}
((da)f)(m)= (d(af))(m)&  =  & d(af(m)) - af(dm)   \\
          &  =  &  (\partial a)f(m) + a d(f((m)) - af(dm)  \\
          &  =  &  (a(df))(m) + (\partial a)f(m)  \\
          &  =  &  ((ad + \partial a)f)(m)
\end{array}     \] 
We finally get $(d_if)^k_{\mu}=(d_if)(y^k_{\mu})={\partial}_if^k_{\mu}-f^k_{\mu +1_i}$ and thus recover 
{\it exactly} the Spencer operator of the previous section though {\it this is not evident at all}. We also get 
$(d_id_jf)^k_{\mu}={\partial}_{ij}f^k_{\mu}-{\partial}_if^k_{\mu+1_j}-{\partial}_jf^k_{\mu+1_i}+f^k_{\mu+1_i+1_j} \Longrightarrow d_id_j=d_jd_i, \forall i,j=1,...,n$ and thus $d_iR_{q+1}\subseteq R_q\Longrightarrow d_iR\subset R$ induces a well defined operator $R\longrightarrow T^*\otimes R:f \longrightarrow dx^i\otimes d_if$. This operator has been first introduced, up to sign, by F.S. Macaulay as early as in $1916$ but this is still not acknowledged ([12]). For more details on the Spencer operator and its applications, the reader may look at ([19],[24],[30],[29]).  \\
\hspace*{13cm}  $\Box$ \\

\noindent
{\bf PROPOSITION 2.2.11}: When $M$ and $N$ are left $D$-modules, then $M {\otimes}_KN$ is also a left $D$-module.  \\

\noindent
{\it Proof}: As before, we may define:  \\
\[    a(m\otimes n)=am \otimes n = m \otimes an, \hspace{5mm}  \forall a \in K, \forall m \in M , \forall n \in N  \]
\[    \xi (m\otimes n)=\xi m \otimes n + m \otimes \xi n, \hspace{5mm}  \forall \xi \in T, \forall m \in M , \forall n \in N   \]
and let the reader finish as an exercise.  \\
\hspace*{13cm}   $\Box$  \\

\noindent
{\bf COROLLARY 2.2.12}: The two structures of left $D$-modules obtained in these two propositions are coherent with the following {\it adjoint isomorphism} existing for any triple $ L, M, N \in mod(D)$:  \\
\[    hom_D(M{\otimes}_K N , L) \stackrel{\varphi}{\longrightarrow} hom_D(M, hom_K(N,L))   \]

\noindent
{\it Proof}: Whenever $f\in hom_D(M {\otimes}_KN,L)$, we may {\it define} ${\varphi}(f)=g$ by $(g(m))(n)=f(m\otimes n) \in L$ and we have successively for any 
$\xi \in T$:  
\[  (\xi(g(m)))(n) = \xi((g(m))(n)) - (g(m))(\xi n)=\xi(f(m\otimes n)) -f(m\otimes \xi n)=f(\xi(m\otimes n)) - f(m \otimes \xi n)  \]
that is $(\xi(g(m)))(n) = f(\xi m \otimes n) =(g(\xi m))(n) $ and thus $\xi (g(m))=g(\xi m), \forall m \in M$. \\
The inverse morphism can be studied similarly.  \\
\hspace*{13cm}  $\Box$  \\

 \noindent
 {\bf COROLLARY 2.2.13 }: $R= hom_K(M,K)  \simeq hom_D(M, hom_K(D,K))$.  \\
 
 \noindent
 {\it Proof}: As $K$ is a field, thus a commutative ring, we have the isomorphism of left $D$-modules $M{\otimes}_K N \simeq N {\otimes}_K M: m \otimes n \rightarrow n \otimes m, \forall m \in M, \forall n \in N$ and we may exchange $M$ and $N$. As $K$ is a left differential module for the rule $(P, a) \rightarrow P(a), \forall P \in D, \forall a \in K$, we obtain:   \\
 \[  \begin{array}{rcl}
  hom_D(M , (hom_K(D,K))\simeq hom_D(M {\otimes}_KD, K) &\simeq & hom_D(D {\otimes }_K M, K) \\
  &  \simeq &  hom_D(D,hom_K(M,K))  \\
   &   \simeq   &  hom_K(M,K) 
  \end{array}   \]
\hspace*{13cm}   $\Box$  \\ 
 
 \noindent
 {\bf COROLLARY  22.14}: The differential module $hom_K(D,K)$ is an injective differential module.    \\
 
 \noindent
{\it Proof}: When $ 0 \rightarrow M' \rightarrow M \rightarrow M" \rightarrow 0$ is a short exact sequence of modules and $N$ is any module, we have {\it only} the exact sequence $ 0 \rightarrow hom(M",N) \rightarrow  hom(M, N) \rightarrow hom(M',N) $ obtained by composition of morphisms. In the present situation, using the previous corollaries, we have the following commutative and exact diagram because $K$ is a field (See [4], p 18):  \\
\[  \begin{array}{cccl}
  0  &    & 0 &    \\
  \downarrow &  &  \downarrow  &   \\
  hom_D(M, hom_K(D,K)) &  \rightarrow &  hom_D(M', hom_K(D,K)) &    \\
     \downarrow &  &  \downarrow  &   \\
   hom_K(M,K) & \rightarrow  & hom_K(M', K)  &  \rightarrow 0  \\
   \downarrow &   &  \downarrow  &   \\
   0 &  &  0 &   
   \end{array}  \]
Chasing in this diagram, we deduce that the upper morphism is an epimorphism and $hom_K(D,K)$ is an injective module because $hom(\bullet, hom_K(D,K))$ transforms a short exact sequence into a short exact sequence. The reader may compare such an approach with the one used in ([41] or [42]) in order to understand why these applications are not dealing with variable coefficients as the differential structure on $R$ {\it must} be defined by the Spencer operator but we do not know any other reference (compare to [14]).   \\
\hspace*{13cm}    $\Box$  \\

\noindent 
{\bf DEFINITION 2.2.15}: With any differential module $M$ we shall associate the {\it graded module} $G=gr(M)$ over the polynomial ring $gr(D)\simeq K[\chi]$ by setting $G={\oplus}^{\infty}_{q=0} G_q$ with $G_q=M_q/M_{q+1}$ and we get $g_q=G_q^*$ where the {\it symbol} $g_q$ is defined by the short exact sequences: \\
\[ 0\longrightarrow M_{q-1}\longrightarrow M_q \longrightarrow G_q \longrightarrow 0  \hspace{4mm}  \Longleftrightarrow \hspace{4mm}  0 \longrightarrow g_q \longrightarrow R_q \longrightarrow R_{q-1} \longrightarrow 0  \]
We have the short exact sequences $0\longrightarrow D_{q-1} \longrightarrow D_q \longrightarrow S_qT \longrightarrow 0 $ leading to $gr_q(D)\simeq S_qT$ and we may set as usual $T^*=hom_K(T,K)$ in a coherent way with differential geometry. \\

The two following definitions, which are well known in commutative algebra, are also valid (with more work) in the case of differential modules (See [21] for more details or the references ([13],[22],[39],[40]) for an introduction to homological algebra and diagram chasing).  \\

\noindent
{\bf DEFINITION 2.2.16}: The set of elements $t(M) = \{ m \in M \mid \exists 0 \neq P \in D, Pm=0\}\subseteq M$ is a differential module called the {\it torsion submodule} of $M$. More generally, a module $M$ is called a {\it torsion module} if $t(M)=M$ and a {\it torsion-free module} if $t(M)=0$. In the short exact sequence $0 \rightarrow t(M) \rightarrow M \rightarrow M' \rightarrow 0$, the module $M'$ is torsion-free. Its defining module of equations $I'$ is obtained by adding to $I$ a representative basis of $t(M)$ set up to zero and we have thus $I \subseteq I'$.  \\

\noindent
{\bf DEFINITION 2.2.17}: A differential module $F$ is said to be {\it free} if $F \simeq D^r$ for some integer 
$r > 0$ and we shall {\it define} $rk_D(F)=r$. If $F$ is the biggest free differential module contained in $M$, then $M/F$ is a torsion differential module and $hom_D(M/F,D)=0$. In that case, we shall {\it define} the {\it differential rank} of $M$ to be $rk_D(M)=rk_D(F)=r$. Accordingly, if $M$ is defined by a linear involutive operator of order $q$, then $rk_D(M)={\alpha}^n_q$.   \\

\noindent
{\bf PROPOSITION 2.2.18}: If $0 \rightarrow M' \rightarrow M \rightarrow M" \rightarrow 0$ is a short exact sequence of differential modules and maps or operators, we have $rk_D(M)=rk_D(M') + rk_D(M")$. \\

\noindent
{\bf REMARK  2.2.19}: We emphasize once more that the left $D$-module $hom_K(D,K)$ used in the literature ([12,[14]°541],[42],[45]) is coming from the right action of $D$ on any $f\in hom_K(D,K)$ through the formula $(Qf)(P)=f(PQ)$ and we have thus:  \\
\[  (af)(P)=f(Pa), \hspace{1cm}  (\xi f)(P)= f(P \xi), \hspace{1cm} \forall a \in K, \forall \xi \in T, \forall P \in D  \]
As $ hom_K(M,K)\simeq hom_D(M,hom_K(D,K))$, they are not coherent at all with the formulas of Proposition 2.2.1, namely:   \\
\[     (af)(m)=a(f(m))=f(am), \hspace{1cm}  (\xi f)(m)= \xi (f(m)) - f(\xi m), \hspace{1cm} \forall a \in K, \forall \xi \in T, \forall m \in M  \]
{\it unless} $K$ {\it is a field of constants}, in particular because, when $M=D$, then $Pa\neq a P$ in general ([30], p 66 for more details). Accordingly, most of the applications of differential duality to control theory must be therefore revisited with these new methods of differential homological algebra (Compare to [45]). We also claim that the use of the adjoint operator must become essential for all the applications to mathematical physics.  \\.

In order to conclude this section, we may say that the main difficulty met when passing from the differential framework to the algebraic framework is the " {\it inversion} " of arrows. Indeed, with $dim(E)=m, dim(F)=p$, when an operator ${\cal{D}}:E \rightarrow F$ is injective, that is when we have the exact sequence $ 0 \rightarrow E \stackrel{{\cal{D}}}{\longrightarrow} F$, like in the case of the operator $0 \rightarrow E \stackrel{j_q}{\longrightarrow} J_q(E) $, on the contrary, using differenial modules, we have the epimorphism $D^p \stackrel{{\cal{D}}}{\longrightarrow} D^m \rightarrow 0$. The case of a formally surjective operator, like the $div$ operator, described by the exact sequence $E \stackrel{{\cal{D}}}{\longrightarrow} F \rightarrow 0$ is now providing the exact sequence of differential modules $ 0 \rightarrow D^p \stackrel{{\cal{D}}}{\longrightarrow} D^m \rightarrow M  \rightarrow 0$ because ${\cal{D}}$ has no CC.  \\

We are now ready for using the results of the second section on the Cartan-K\"{a}hler theorem. For such a purpose, separating the {\it parametric} jets that can be chosen arbitrarily from the {\it principal} jets that can be obtained by using the fact that the given OD or PD equations have coefficients in a differential field $K$, we may write the solved equations in the symbolic form $y_{pri} - c^{par}_{pri} y_{par}=0$ with $c\in K$ and an implicit (finite) summation in order to obtain for the sections $f_{pri} - c^{par}_{pri}f_{par}=0$. Using the language of Macaulay, it follows that the so-called {\it modular equations} are $E\equiv f_{pri}a^{pri} + f_{par}a^{par}=0$ with eventually an infinite number of terms in the implicit summations. Substituting, we get at once $E\equiv f_{par}(a^{par} + c^{par}_{pri}a^{pri})=0$. Ordering the $y_{par}$ as we already did and using a basis $\{ (1,0,...), (0,1,0,...), (0,0,1,0,...), ... \}$ for the $f_{par}$, we may select the {\it parametric modular equations} $E^{par}\equiv a^{par} + c^{par}_{pri}a^{pri}=0$. \\

When $k$ is a field of constants, a polynomial $P=a^{\mu}{\chi}_{\mu}\in k[\chi]$ of degree $q$ is multiplied by a monomial ${\chi}_{\nu}$ with $\mid \nu \mid=r$, we get ${\chi}_{\nu}P=a^{\mu}{\chi}_{\mu+\nu}$. Hence, if $0\leq \mid \mu\mid \leq q$, the "{\it shifted} " polynomial thus obtained is such that $r\leq \mid \mu+\nu\mid \leq q+r$ and the difference between the maximum degree and the minimum degree of the monomials involved is always equal to $q$ and thus fixed. When $n=1$, one can exhibit a series only made with $0$ or $1$ like $f=(1,0,1,0,0,1,0,0,0,1,0,0, ...)$ with " {\it zero zones} " of successive increasing lengths $1,2,3,4, ...$ and so on, separated by $1$ in such a way that the contraction with the shifted polynomial is the leading term of the given polynomial and extend this procedure to $n$ arbitrary.  \\ 

Replacing ${\chi}_i$ by $d_i$ and {\it degree} by {\it order}, we may use the results of section 2 in order to split the CK-data into $m$ formal power series of $0$ (constants), $1,...,n$ variables that we shall call series of {\it type} $i$ for $i=0,1,...,n$. However, as the following elementary example will show, the shifting procedure cannot be applied to the variable coefficient case, namely when $K$ is used in place of $k$. Indeed, with $n=1, d_x=d, K=\mathbb{Q}(x)$ and $P=d^2-\frac{x}{3}$, if we contract $d^3P=d^5-\frac{x}{3}d^3-d^2$ with the series $f=(1,0,1,0,0,1,0,...)$ already defined when $n=1$, we get $1-1=0$ though $P$ does not kill $f$ because $d^2f=(1,0,0,1,0,0,...)\Rightarrow Pf=(1-\frac{x}{3},...)$ and the contraction of $P$ with $f$ is $1-\frac{x}{3}\neq 0$. \\

WE SHALL ESCAPE FROM THIS DIFFICULTY BY MEANS OF A TRICK BASED ON A SYSTEMATIC USE OF THE SPENCER OPERATOR (Compre to [14]).  \\

The idea will be to shift the series to the {\it left} ({\it decreasing ordering}), up to sign, instead of shifting the operator to the {\it right} ({\it increasing ordering}). For this, we notice that we want that the contraction of $P=a^{\mu}d_{\mu}$ where $\mid \mu \mid \leq q=ord(P)$ with $f$ should be zero, that is $a^{\mu}f_{\mu}=0\Rightarrow ({\partial}_ia^{\mu})f_{\mu}+a^{\mu}({\partial}_if_{\mu})=0,\forall i=1,...,n$. But $d_iP=a^{\mu}d_{\mu+1_i}+({\partial}_ia^{\mu})d_{\mu}$ must also contract to zero wih $f$ that is $a^{\mu}f_{\mu+1_i}+({\partial}_ia^{\mu})f_{\mu}=0$. Substracting, we obtain therefore the condition $a^{\mu}({\partial}_if_{\mu}-f_{\mu+1_i})=0$, that is $P$ must also contract to zero with the shift $d_if$ or even $d_{\nu}f$ of $f$ when $f$ is made with $0$ and $1$ only. Applying this computation to the above example, we get $-df=(0,1,0,0,1,0,...)\Rightarrow d^2f=(1,0,0,1,0,...)\Rightarrow -d^3f=(0,0,1,0,...)$ and the contraction with $P$ provides the leading coefficient $1\neq 0$ of $P$ like the contraction of $d^3P$ with $d^4f = (0,1,0,0,0,1,0,...) $, that is the same series can be used but in a quite different framework. Also, in the finite dimensional case existing when the symbol $g_q$ of $R_q$ is finite type, that is when $g_{q+r}=0$ for a certain integer $r\geq 0$, applying the $\delta$-sequence inductively to $g_{q+n+i}$ for $i=r-1,..., 0$ as in ([17], Proposition 4.7, p 123), it is known that $g_q$ is finite type {\it and} involutive if and only if $g_q=0$, that is to say $dim(R_q)=dim(R_{q-1})$. In a coherent way, we have thus obtained:   ÊÊ\\

\noindent
{\bf THEOREM 2.2.20}: If $M$ is a differential module over $D=K[d]$ defined by a first order involutive system in the $m$ unknowns $y^1,...,y^m$ with no zero order equation, the differential module $R=hom_K(M,K)$ may be generated over $D$ by a {\it finite basis of sections} containing $m$ generators. \\

In the general situation, counting the number of CK data, we have ${\alpha}^1_q+...+{\alpha}^n_q=dim(g_q)$ and $dim(R_q)=dim(g_q)+dim(R_{q-1})$. We obtain therefore the following result which is coherent with the number of unknowns in the Spencer form $R_{q+1}\subset J_1(R_q)$.  \\

\noindent
{\bf COROLLARY 2.2.21}: If $M$ is a differential module over $D=K[d]$ defined by an involutive system $R_q\subset J_q(E)$, the differential module $R=hom_K(M,K)$ may be generated over $D$ by a {\it finite basis of sections} containing $dim(R_q)$ generators.  \\

\noindent
{\bf EXAMPLE  2.2.22}: Among the most interesting examples with $m=1, n=3, K= \mathbb{Q}$ we present the following second order system provided by Macaulay in $1916$ ([ M], p   ):  \\
\[         y_{33} =0, \hspace{2cm}   y_{23} - y_{11} =0 , \hspace{2cm}   y_{22}=0   \]
It is easy to check that $g_2$ with $dim(g_2)=3$ is not involutive, that $g_3$ with $dim(g_3)=1$ is $2$-acyclic because $y_{123} - y_{111}=0$ and that $g_4=0$ is trivially involutive. Accordingly, only the system $R_4$ is involutive with the $8=2^3$ parametric jets $(y, y_1,y_2,y_3,y_{11},y_{12}, y_{13}, y_{111})$ and all the three characters vanish both with all the jets of order $\geq 4$ because the system is homogeneous.  \\
Let us consider the following $8$ sections with all other components equal to zero:  \\
\noindent
\[  \begin{array}{|c|c|c|c|c|c|c|c|c|c|c|}
  section& y & y_1 & y_2 & y_3 & y_{11} & y_{12} & y_{13} & y_{23} & y_{111} & y_{123}   \\
\hline
f_1& 1 & 0 & 0 & 0 & 0 & 0 & 0 & 0 & 0 & 0  \\
f_2 & 0 & 1 & 0 & 0 & 0 & 0 & 0 & 0  &0 & 0 \\
f_3 & 0 & 0 & 1 & 0 & 0 & 0 & 0 & 0 & 0 & 0  \\
f_4  & 0& 0 & 0 & 1 & 0 & 0 & 0  & 0  &0  & 0  \\
 f_5 & 0 & 0 & 0 &  0  & 1 & 0 & 0  & 1 & 0 & 0  \\
 f_6& 0 & 0 & 0 & 0 & 0 & 1 & 0 & 0 & 0 & 0  \\
 f_7& 0 & 0 & 0 & 0 &  0 & 0 & 1 & 0 & 0 & 0  \\
 f_8& 0 & 0 & 0 & 0 & 0 & 0 & 0 & 0 &1 &  1  
\end{array}  \]      
Taking into account the two PD equations $ y_{23} - y_{11}=0$ and $y_{123} - y_{111}=0$, we obtain successively:  \\
\[     d_1f_8= - f_5, \,\,d_1f_7= - f_4, \,\, d_3f_7= - f_2, ............, d_1f_2= - f _1     \]
Or, equivalently, working with the corresponding modular equations:  \\
\[   E_8 \equiv a^{111} + a^{123}=0  \,\,\,  \Rightarrow  \,\,\, - d_1E_8 \equiv E_5 \equiv a^{11} + a^{23}=0  \]
and so on. It follows that all the sections can be generated by the single section $f_8$ and all the modular equations can be generated by the single modular equation $E_8=0$, a result absolutely not evidet at first sight but coherent with the fact that the radical of the annihilator of $M$ is the maximal ideal 
$\mathfrak{m}=(d_1,d_2,d_3)$.  \\
Finally, with $m=1, n=3$ but $K=\mathbb{Q}(x^1,x^2,x^3)$ an even more striking example has been provided by M. Janet in $1920$ ([6]) with the following second order system:  \\
\[    y_{33} - x^2 \, y_{11} = 0 , \hspace{2cm}  y_{22}=0  \]
In this case, $dim(R)=12 < \infty $ but $R$ can be generated by the single modular equation:    \\
\[        E \equiv a^{12333} + x^2 \, a^{1333} + a^{1113} = 0   \]
because all the jets of order $> 5$ vanish (See [20] and [22] for more details).    \\

\noindent
{\bf EXAMPLE 2.2.23}: If $K=\mathbb{Q}(x)$ and $n=1,m=1$, let us consider the third order OD equation $\Phi \equiv y_{xxx}-y_x=0$ for which we may exhibit the basis of sections:     \\
\[  \{ f_1=(1,0,0,0,0, ...),f_2=(0,1,0,1,0,1,...),f_3=(0,0,1,0,1,0,...)\}    \]. 
With $d=d_x$ and $\partial = {\partial}_x$, we obtain $ df_1=0, df_2= - f_1 - f_3, df_3= - f_2 $ and check that all the sections can be generated by a single one, namely $f_3$ which describes the power series of the solution $y=f(x)=ch(x)-1$. We have indeed $\partial f=sh(x)$ and ${\partial}^2f(x)=ch(x)$.   \\
With now $m=2$, let us consider the differential module defined by the system $y^1_{xx} - y^1=0,y^2_x=0$. Setting $y=y^1-y^2$, we successively get:   \\
\[  \begin{array}{rcl}
y_x=y^1_x, y_{xx}=y^1_{xx}=y^1, y_{xxx} = y^1_x  & \Rightarrow & 
y^1=y_{xx}, y^1_x=y_x, y^1_{xx}=y_{xx}, y^1_{xxx}=y_x, ... \\
 & \Rightarrow  &  y^2=y_{xx}-y, y^2_x=0, ...
 \end{array}    \]
and a differential isomorphism with the module defined by the new system $y_{xxx}-y_x=0$. We have seen that the sections of the second system are easily seen to be generated by the single section $f_3=(0,0,1,0, ...)$, a result leading to the only generating section $f^1= 1, f^1_x= 0, f^1_{xx}=1, ..., f^2=1, f^2_x=0,...$ of the initial system but {\it these sections do not describe solutions} because ${\partial}f^1-f^1_x=  0$ and ${\partial}f^2-f^2_x=0$ but ${\partial}f^1_x -f^1_{xx}= - 1 \neq 0$. {\it We do not know any reference in computer algebra dealing with sections} (See [26] for more details)  \\

\noindent
{\bf EXAMPLE 2.2.24}: With $n=1, m=1, q=2, K=\mathbb{Q}(x)$, let us consider the second order OD equation  $\Phi \equiv y_{xx}-xy=0$. We successively obtain by prolongation:  \\
\[   y_{xxx}-xy_x-y=0, y_{xxxx}-2y_x-x^2y=0, y_{xxxxx}-x^2y_x-4xy=0, y_{xxxxxx}-6xy_x-(x^3+4)y=0 \] 
and so on. We obtain the corresponding board describing the maps ${\rho}_r(\Phi)$ (Compare to [M]):Ê \\
\[  \begin{array}{r|c|c|c|c|c|c|c|l}
 order & y & y_x & y_{xx} & y_{xxx} & y_{xxxx} & y_{xxxxx} & y_{xxxxxx} &... \\
\hline
2 & -x & 0 & 1 & 0 & 0 & 0 & 0& ...   \\
3 & -1 & -x & 0 & 1 & 0 & 0& 0& ...   \\
4& -x^2 & -2  & 0 & 0 & 1 & 0& 0& ... \\
5& -4x & -x^2 & 0&0 & 0 & 1 & 0 & ... \\
6 & -(x^3+4) & -6x & 0 & 0 & 0 & 0 & 1 & ...
\end{array}  \]

Let us define the sections $f'$ and $f"$ by the following board where $d=d_x$: \\
\[  \begin{array}{r|c|c|c|c|c|c|c|l}
  section& y & y_x & y_{xx} & y_{xxx} & y_{xxxx} & y_{xxxxx} & y_{xxxxxx} &... \\
\hline
f' & 1 & 0 & x & 1 & x^2 & 4x & x^3+4& ...   \\
f" & 0 & 1 & 0 & x & 2 & x^2 & 6x& ...   \\
\hline
df' & 0 & -x  & 0 & -x^2 & -2x & - x^3 & -6x^2 & ... \\
df" & -1 & 0 & -x & -1 & -x^2 & -4x & -x^3-4 & ... \\
\end{array}  \]
in order to obtain $df'=-xf", df"=-f'$. Though this is not evident at first sight, {\it the two boards are orthogonal over} $K$ in the sense that each row of one board contracts to zero with each row of the other though only the rows of the first board do contain a finite number of nonzero elements. {\it It is absolutely essential to notice that the sections} $f'$ {\it and} $f"$ {\it have nothing to do with solutions}  because $df'\neq 0, df"\neq 0$ on one side and also because $d^2f'-xf'=-f"=\frac{1}{x}df'\neq 0$ even though $d^2f"-xf"=0$ on the other side. As a byproduct, $f'$ or $f"$ can be chosen separately as unique generating section of the inverse system over $K$ ({\it care}) and we may write for example:   \\
\[ f' \rightarrow E'\equiv a^0+xa^{xx}+a^{xxx}+x^2 a^{xxxx}+ ... =0, \,\,\, f" \rightarrow E"\equiv a^x+xa^{xxx}+2a^{xxxx}+ ... =0 \]     \\

\noindent
{\bf EXAMPLE 2.2.25}: With $n=1, m=2, q= 2, K=\mathbb{Q}(x)$, let us consider again the second order system $y^1_{xx} - y^1=0,y^2_x=0$. Setting  $z^1=y^1, z^2=y^1_x, z^3=y^2$, we obtain the first order involutive system:
\[  \left\{   \begin{array}{lcc}
z^1_x-z^2 & = & 0 \\
z^2_x - z^1& =  & 0 \\
z^3_x  & = & 0 
\end{array}
\right. \fbox{$\begin{array}{c}
 x \\
 x  \\
 x
 \end{array}$}  \]
It follows that the CK data for $z=g(x)$ are $\{g_1=g^1(0), g_2=g^2(0), g_3=g^3(0)\}$. Using the given equations and their solved prolongations like $y^1_{xx} -y^1=0, y^2_{xx} - y^2=0$ and so on, we have the finite basis ({\it care} !): \\
\[  \begin{array}{l|ccc|ccc|ccc|c}
    & z^1 & z^2 & z^3 & z^1_x & z^2_x & z^3_x & z^1_{xx} &z^2_{xx} & z^3_{xx} & ... \\
\hline
g_1  & 1 & 0 & 0 & 0 & 1 & 0 & 1 & 0 & 0 & ... \\
g_2 & 0 & 1 & 0 & 1 & 0 & 0 & 0 & 1 & 0 & ... \\
g_3 & 0 & 0 & 1 & 0 &0 & 0 & 0 & 0 & 0 & ... 
\end{array}  \] 
As $dg_1= - g_2,dg_2= - g_1, dg_3=0$, a basis with only two generators may be $\{g^2,g^3\}$. However: \\
\[   h=g_1- g_3, dh= - g_2, d^2h=g_1  \Leftrightarrow g_1=d^2h, g_2= - dh, g_3=d^2h-h  \]
and we obtain the unique generator $h$ (See [26] for details).     \\

\noindent
{\bf 2.3) LINEAR CONTROL THEORY}  \\

The most striking aspect of the application of module/system theory to  linear control theory is that it is coming from rather unexpected chases in commutative and exact diagrams looking like rather abstract at first sight. As more details and examples can be found in book form ([21],[22]), we shall only provide below a few new results that cannot be found elsewhere.  \\

 \noindent
{\bf PROPOSITION 2.3.1}: One has the short exact sequence  of (differential) modules:   \\
\[       0 \longrightarrow  (L + M) \longrightarrow N  \longrightarrow  (N/L)/(M/(L\cap M)) \longrightarrow 0   \]

\noindent
{\it Proof}: Using elementary classical homological algebra, one obtains the following commutative and exact diagram:   \\

\[ \begin{array}{rcccccl}
   & 0  & &  0  & &  0 &  \\
   &  \downarrow &  &  \downarrow  &  &  \downarrow  &   \\
0 \rightarrow & L\cap M & \rightarrow & L & \rightarrow & L/(L\cap M) & \rightarrow 0  \\
  & \downarrow&  \searrow  & \downarrow &  &  \downarrow  \\
0 \rightarrow & M & \rightarrow & N & \rightarrow & N/M & \rightarrow 0   \\
    &  \downarrow &  &  \downarrow  & \searrow & \downarrow &    \\
0 \rightarrow  & M / (L\cap M) & \rightarrow & N/L & \rightarrow & (N/L)/(M/(L\cap M)) &  \rightarrow 0   \\
  &  \downarrow  &  &  \downarrow  &  & \downarrow  \\
   &   0  &  &  0  &  &  0  &  
\end{array}  \]
 At first, the lower southeast arrow being the composition of two epimorphisms is an epimorphism. A circular chase finally proves that any element of $N$ killed by this southeast arrow is the sum of an element of $L$ and an element of $M$, achieving the proof. It is important to notice the symmetric part plaid by $L$ and $M$ in $N$.  \\
\hspace*{12cm}  $\Box$  \\
 In the general situation, we obtain from the left upper commutative square the useful formulas:  \\
 \[   rk_D(N/(L\cap M)=rk_D(L/(L\cap M)) + rk_D(N/L)=rk_D(M/(L\cap M)) + rk_D(N/M)   \]
 in a coherent way with the following corollary:  \\
 
 \noindent
 {\bf COROLLARY 2.3.2}: If $L+M=N$, then one has $L/(L\cap M) \simeq N/M$ and $M/(L\cap M)\simeq N/L$.  \\
 
 \noindent
 {\bf THEOREM 2.3.3}: There is a bijective correspondence between the intermediate differential modules $ (L\cap M) \subset L' \subset L$ and the intermediate differential modules 
 $M \subset N' \subset N$ defined by the rules:  \\
 \[   L' \rightarrow (L'+ M )= N' \subset N, \,\,\,\,\,    N' \rightarrow (L \cap N') = L' \subset L   \]
 
\noindent
{\it Proof}: We have the following commutative diagram of injections:   \\

\[   \begin{array}{ccccl}
M & \rightarrow &  N' & \rightarrow &  N  \\
\uparrow &   &   \uparrow  &  &  \uparrow  \\
L\cap M & \rightarrow &  L'  & \rightarrow  & L 
\end{array}  \]
Let us start with $L'$, construct $N'= L'+M$ and obtain $L"=L\cap N'$. First of all, we get $L'\subset L$ and $L' \subset N'  \Rightarrow L'\subseteq L"$. Now, using the left commutative square, we obtain from the previous proposition $N' / L'\simeq M / (L\cap M)$. Similarly, using the right commutative square, we obtain $N' / N" \simeq N / L$ and thus an isomorphism $N' / L' \simeq N' / L"$. However, we have the following commutative and exact diagram:  \\
\[ \begin{array}{rcccccl}
 & 0  & &  0  & &  0 &  \\
   &  \downarrow &  &  \downarrow  &  &  \downarrow  &   \\
0 \rightarrow & L' & = & L' & \rightarrow & 0 \\
  & \downarrow &  & \downarrow &  &  \downarrow  \\
0 \rightarrow & L" & \rightarrow & N' & \rightarrow & N' / L" & \rightarrow 0   \\
    &  \downarrow &  &  \downarrow  &   & \parallel &    \\
0 \rightarrow  & L" / L' & \rightarrow & N' / L' & \rightarrow & N' / L" &  \rightarrow 0   \\
  &  \downarrow  &  &  \downarrow  &  & \downarrow  \\
   &   0  &  &  0  &  &  0  &  
\end{array}  \]
and thus $L" / L' = 0$ that is $L' = L"$.  \\  
Finally, starting with $N'$, we should obtain $L'=N' \cap L $, then define $N"= L' + M \subseteq N'$ and conclude as before that $N' = N"$. Replacing specialzations by injections while chasing in the following commutative and exact diagram:  \\
\[ \begin{array}{rcccccl}
0 \rightarrow &  N' &  \rightarrow & N  &  \rightarrow &  N"  &    \rightarrow  0 \\
                   &    \uparrow  &  &  \uparrow  &  &  \uparrow  &        \\
0 \rightarrow &  L' &  \rightarrow & L &  \rightarrow &  L"  &    \rightarrow   0\\
                  &    \uparrow  &  &  \uparrow  &  &  \uparrow  &        \\
                  & 0    &  &  0  &  &  0  &
\end{array}   \]
we have thus been able to deal only with submodules of $N$. In paricular, if $N$ is torsion-free, that is $t(N)=0$, the interest of this aproach is that all the submodules are torsion-free. \\

\hspace*{12cm}   $\Box$   \\ 
 
\noindent
{\bf EXAMPLE 2.3.4}: (See [21], p 736-738] for the details and diagrams) With $n=1, K= \mathbb{Q} $, let us consider the differential module defined by the OD equation $d^5v +3 d^4 v - 4 d^2v=d^3u- du$. We may define the input differential module $M_{in}$ by using $u$ and the output differential module by using 
$y=d^2 v$. The differential module $(M_{in} + M_{out}) \subset M$ with a strict inclusion, is defined by the OD equation $d^3y + 3d^2y - 4 y=d^3u- du$. Localizing at the covector $\chi$, we get the equation $(\chi - 1)((\chi) + 2)^2y = \chi (\chi -1 )(\chi + 1)u$ that we can also write $(d-1)(d+2)^2v=d(d-1)(d+1)u$ because $K$ is a field of constants. As we can factor by $(d-1)$ it follows that $t(M_{in} + M_{out})$ is generated by $z= d^2 y + 4 dy + 4y -d^2 u + du$ that satisfies $dz + z=0$. We have $ann(M/M_{in})=(d^2 (d-1)(d+2)^2)$ and its radical is $rad(M/M_{in})=(d(d_1)(d+2))=(d) \cap (d-1) \cap (d+2)$ is an intersection of prime ideals. Similarly, we have $y=0 \Rightarrow d^2v=0, (d^3 -d)u=0$ and thus $ann(M/M_{out})= (d^2) \cap (d^3 - d)$ leading to $rad(ann(M/M_{out}))= (d) \cap (d-1) \cap (d+1)$ as an intersection of prime differential ideals. We have proved in ([21]) how to use these differential submodules of $M$ both with the new differential modules $M'_{in}=t(M) + M_{in}$ and $M'_{out}=t(M) + M_{out}$ in order to study all the problems concerning poles and zeros. In the present paper, as we are only interested by controllability, we have just to study the differential submodules of the torsion-free differential module $M/t(M)$.\\

\noindent
{\bf 3) NONLINEAR CORRESPONDENCES}  \\

\noindent
{\bf 3.1) NONLINEAR SYSTEMS}  \\

If $X$ is a manifold with local coordinates $(x^i)$ for $i=1, ... ,n=dim(X)$, let us consider the {\it fibered manifold} ${\cal{E}}$ over $X$ with ${dim}_X({\cal{E}})=m$, that is a manifold with local coordinates $(x^i,y^k)$ for $i=1,...,n$ and $k=1,...,m$ simply denoted by $(x,y)$, {\it projection} $\pi:{\cal{E}}\rightarrow X:(x,y)\rightarrow (x)$ and changes of local coordinates $\bar{x}=\varphi(x), \bar{y}=\psi(x,y)$. If $\cal{E}$ and $\cal{F}$ are two fibered manifolds over $X$ with respective local coordinates $(x,y)$ and $(x,z)$, we denote by ${\cal{E}}{\times}_X{\cal{F}}$ the {\it fibered product} of $\cal{E}$ and $\cal{F}$ over $X$ as the new fibered manifold over $X$ with local coordinates $(x,y,z)$. We denote by $f:X\rightarrow {\cal{E}}: (x)\rightarrow (x,y=f(x))$ a global {\it section} of $\cal{E}$, that is a map such that $\pi\circ f=id_X$ but local sections over an open set $U\subset X$ may also be considered when needed. Under a change of coordinates, a section transforms like $\bar{f}(\varphi(x))=\psi(x,f(x))$ and, differentiating with respect to $x^i$, we may introduce new coordinates $(x^i,y^k,y^k_i)$ transforming like:\\
\[ {\bar{y}}^l_r{\partial}_i{\varphi}^r(x)=\frac{\partial{\psi}^l}{\partial x^i}(x,y)+\frac{\partial {\psi}^l}{\partial y^k}(x,y)y^k_i  \]
We shall denote by $J_q({\cal{E}})$ the {\it q-jet bundle} of $\cal{E}$ with local coordinates $(x^i, y^k, y^k_i, y^k_{ij},...)=(x,y_q)$ called {\it jet coordinates} and sections $f_q:(x)\rightarrow (x,f^k(x), f^k_i(x), f^k_{ij}(x), ...)=(x,f_q(x))$ transforming like the sections $j_q(f):(x) \rightarrow (x,f^k(x), {\partial}_if^k(x), {\partial}_{ij}f^k(x), ...)=(x,j_q(f)(x))$ where both $f_q$ and $j_q(f)$ are over the section $f$ of $\cal{E}$. It will be useful to introduce a {\it multi-index} $\mu=({\mu}_1, ... ,{\mu}_n)$ with length $\mid \mu \mid={\mu}_1+ ... +{\mu}_n$ and to set ${\mu}+1_i=({\mu}_1  ... ,{\mu}_{i-1}, {\mu}_i+1, {\mu}_{i+1},...,{\mu}_n)$. Also, a jet coordinate $y^k_{\mu}$ is said to be of {\it class} $i$ if ${\mu}_1=...={\mu}_{i-1}=0, {\mu}_i\neq 0$. As the background will always be clear enough, we shall use the same notation for a vector bundle or a fibered manifold and their sets of sections [31,36]. We finally notice that $J_q({\cal{E}})$ is a fibered manifold over $X$ with projection ${\pi}_q$ while $J_{q+r}({\cal{E}})$ is a fibered manifold over $J_q({\cal{E}})$ with projection ${\pi}^{q+r}_q, \forall r\geq 0$ [, , ].\\

\noindent
{\bf DEFINITION 3.1.1}: A ({\it nonlinear}) {\it system} of order $q$ on $\cal{E}$ is a fibered submanifold ${\cal{R}}_q\subset J_q({\cal{E}})$ and a global or local {\it solution} of ${\cal{R}}_q$ is a section $f$ of $\cal{E}$ over $X$ or $U\subset X$ such that $j_q(f)$ is a section of ${\cal{R}}_q$ over $X$ or $U\subset X$.\\

\noindent
{\bf DEFINITION 3.1.2}: When the changes of coordinates have the linear form $\bar{x}=\varphi(x),\bar{y}= A(x)y$, we say that $\cal{E}$ is a {\it vector bundle} over $X$. Vector bundles will be denoted by capital letters $E,F, \dots $ and will have sections denoted by $\xi,\eta,\dots $. In particular, we shall denote as usual by $T=T(X)$ the {\it tangent bundle} of $X$, by $T^*=T^*(X)$ the {\it cotangent bundle}, by ${\wedge}^rT^*$ the {\it bundle of r-forms} and by $S_qT^*$ the {\it bundle of q-symmetric covariant tensors}. When the changes of coordinates have the form $\bar{x}=\varphi(x),\bar{y}=A(x)y+B(x)$ we say that $\cal{E}$ is an {\it affine bundle} over $X$ and we define the {\it associated vector bundle} $E$ over $X$ by the local coordinates $(x,v)$ changing like $\bar{x}=\varphi(x),\bar{v}=A(x)v$. \\

\noindent
{\bf DEFINITION 3.1.3}: If the tangent bundle $T({\cal{E}})$ has local coordinates $(x,y,u,v)$ changing like ${\bar{u}}^j={\partial}_i{\varphi}^j(x)u^i, {\bar{v}}^l=\frac{\partial {\psi}^l}{\partial x^i}(x,y)u^i+\frac{\partial {\psi}^l}{\partial y^k}(x,y)v^k$, we may introduce the {\it vertical bundle} $V({\cal{E}})\subset T({\cal{E}})$ as a vector bundle over $\cal{E}$ with local coordinates $(x,y,v)$ obtained by setting $u=0$ and changes ${\bar{v}}^l=\frac{\partial {\psi}^l}{\partial y^k}(x,y)v^k$. Of course, when $\cal{E}$ is an affine bundle over $X$ with associated vector bundle $E$ over $X$, we have $V({\cal{E}})={\cal{E}}\times_XE$. We have the short exact sequence of vector bundles over ${\cal{E}}$:  
\[ 0 \rightarrow  V({\cal{E}})  \rightarrow T({\cal{E}}) \stackrel{T(\pi)}{\longrightarrow} {\cal{E}}{\times}_XT \rightarrow   0   \]
Accordingly, in variational calculus, the couple $(f,\delta f)$ made by a section $f$ of ${\cal{E}}$ and its {\it variation} $\delta f$ is nothing else but a section of $V({\cal{E}})$ while $\delta f$ has no reason at all to be " {\it small} ".  \\

For a later use, if $\cal{E}$ is a fibered manifold over $X$ and $f$ is a section of $\cal{E}$, we denote by $f^{-1}(V({\cal{E}}))$ the {\it reciprocal image} of $V({\cal{E}})$ by $f$ as the vector bundle over $X$ obtained when replacing $(x,y,v)$ by $(x,f(x),v) $ in each chart, along with the following commutative diagram:   \\
\[  \begin{array}{ccc}
       E  &  \longrightarrow &  V({\cal{E}})  \\
    \delta f \uparrow   \downarrow \hspace{4mm}&  &  \downarrow  \\
       X  & \begin{array}{c}
        \stackrel{f}{\longrightarrow} \\
        \stackrel{\longleftarrow}{\pi}
        \end{array}     &  {\cal{E}}  
 \end{array}  \]
A similar construction may also be done for any affine bundle over ${\cal{E}}$. When the background is clear enough, with a slight abuse of language, we shall sometimes set $E=V(\cal{E})$ as a vector bundle over $\cal{E}$ and call " {\it vertical machinery} " such a useful systematic notation.  \\
Looking at the transition rules of $J_q(\cal{E})$, we deduce easily the following results: \\

\noindent
{\bf PROPOSITION 3.1.4}: $J_q(\cal{E})$ is an affine bundle over $J_{q-1}(\cal{E})$ modeled on $S_qT^*{\otimes}_{\cal{E}}E$ but we shall not specify the tensor product in general.  \\

\noindent
{\bf PROPOSITION 3.1.5}: There is a canonical isomorphism $V(J_q({\cal{E}})) \simeq J_q(V({\cal{E}}))=J_q(E)$ of vector bundles over $J_q(\cal{E})$ given by setting $v^k_{\mu}=v^k_{,\mu}$ at any order and a short exact sequence: \\
\[  0 \rightarrow S_qT^*\otimes E \rightarrow J_q(E) \stackrel{{\pi}^q_{q-1}}{\longrightarrow} J_{q-1}(E) \rightarrow 0  \]
 of vector bundles over $J_q(\cal{E})$ allowing to establish a link with the formal theory of linear systems.  \\

\noindent
{\bf PROPOSITION 3.1.6}: There is an exact sequence: \\
\[  0 \rightarrow {\cal{E}} \stackrel{j_{q+1}}{\longrightarrow} J_{q+1}({\cal{E}}) \stackrel{d}{\longrightarrow} T^*\otimes J_q(E)  \]
where $df_{q+1}=j_1(f_q) - f_{q+1}$ is over $f_q$ with components $(df_{q+1})^k_{\mu,i}={\partial}_if^k_{\mu}-f^k_{\mu + 1_i}$ is called the (nonlinear) {\it Spencer operator}. \\

\noindent
{\bf DEFINITION 3.1.7}: If ${\cal{R}}_q\subset J_q({\cal{E}})$ is a system of order $q$ on ${\cal{E}}$, then ${\cal{R}}_{q+1}={\rho}_1({\cal{R}}_q)=J_1({\cal{R}}_q)\cap J_{q+1}({\cal{E}})\subset J_1(J_q({\cal{E}}))$ is called the {\it first prolongation} of 
${\cal{R}}_q$ and we may define the subsets ${\cal{R}}_{q+r}$. In actual practice, if the system is defined by PDE ${\Phi}^{\tau}(x,y_q)=0$ the first prolongation is defined by adding the PDE $d_i{\Phi}^{\tau}\equiv {\partial}_i{\Phi}^{\tau} + y^k_{\mu +1_i}{\partial}{\Phi}^{\tau}/\partial y^k_{\mu}=0$. accordingly, $f_q\in {\cal{R}}_q \Leftrightarrow {\Phi}^{\tau}(x,f_q(x))=0$ and $f_{q+1}\in {\cal{R}}_{q+1} \Leftrightarrow {\partial}_i{\Phi}^{\tau}+f^k_{\mu +1_i}(x) \partial {\Phi}^{\tau}/\partial y^k_{\mu}=0$ as identities on $X$ or at least over an open subset $U\subset X$. Differentiating the first relation with respect to $x^i$ and substracting the second, we finally obtain:  \\
\[  ({\partial}_if^k_{\mu}(x) - f^k_{\mu +1_i}(x))\partial {\Phi}^{\tau}/\partial y^k_{\mu}=0 \Rightarrow df_{q+1}\in T^*\otimes R_q \]
and the Spencer operator restricts to $d:{\cal{R}}_{q+1} \rightarrow T^*\otimes R_q$. We set ${\cal{R}}^{(1)}_{q+r}={\pi}^{q+r+1}_{q+r}( {\cal{R}}_{q+r+1})$.  \\

\noindent
{\bf DEFINITION 3.1.8}: The {\it symbol} of ${\cal{R}}_q$ is the family $g_q=R_q \cap S_qT^*\otimes E$ of vector spaces over ${\cal{R}}_q$. The symbol $g_{q+r}$ of ${\cal{R}}_{q+r}$ only depends on $g_q$ by a direct prolongation procedure. We may define the vector bundle $F_0$ over ${\cal{R}}_q$ by the short exact sequence $0 \rightarrow R_q \rightarrow J_q(E) \rightarrow F_0 \rightarrow 0$ and we have the exact induced sequence $0 \rightarrow g_q \rightarrow S_ qT^*\otimes E \rightarrow F_0$ .\\

Setting $a^{\tau \mu}_k(x,y_q)=\partial {\Phi}^{\tau}/\partial y^k_{\mu}(x,y_q) $ whenever $\mid \mu \mid =q$ and $(x,y_q)\in {\cal{R}}_q$, we obtain:  \\
 \[  g_q=\{ v^k_{\mu}\in S_qT^*\otimes E \mid a^{\tau \mu}_k(x,y_q)v^k_{\mu}=0\} , \mid \mu \mid =q, (x,y_q)\in {\cal{R}}_q    \]
  \[  \Rightarrow
 g_{q+r}={\rho}_r(g_q)=\{ v^k_{\mu + \nu}\in S_{q+r} T^*\otimes E \mid a^{\tau \mu}_k(x,y_q)v^k_{\mu + \nu}=0\}, \mid \mu\mid=q, \mid \nu \mid =r, (x,y_q)\in {\cal{R}}_q  \]
In general, neither $g_q$ nor $g_{q+r}$ are vector bundles over ${\cal{R}}_q$.  \\

On ${\wedge}^sT^*$ we may introduce the usual bases $\{dx^I=dx^{i_1}\wedge ... \wedge dx^{i_s}\}$ where we have set 
$I=(i_1< ... <i_s)$. In a purely algebraic setting, one has:  \\

\noindent
{\bf PROPOSITION 3.1.9}: There exists a map $\delta:{\wedge}^sT^*\otimes S_{q+1}T^*\otimes E\rightarrow {\wedge}^{s+1}T^*\otimes S_qT^*\otimes E$ which restricts to $\delta:{\wedge}^sT^*\otimes g_{q+1}\rightarrow {\wedge}^{s+1}T^*\otimes g_q$ and ${\delta}^2=\delta\circ\delta=0$.\\

{\it Proof}: Let us introduce the family of s-forms $\omega=\{ {\omega}^k_{\mu}=v^k_{\mu,I}dx^I \}$ and set $(\delta\omega)^k_{\mu}=dx^i\wedge{\omega}^k_{\mu+1_i}$. We obtain at once $({\delta}^2\omega)^k_{\mu}=dx^i\wedge dx^j\wedge{\omega}^k_{\mu+1_i+1_j}=0$ and $ a^{\tau \mu}_k(\delta \omega)^k_{\mu}=dx^i \wedge(  a^{\tau \mu}_k{\omega}^k_{\mu +1_i})=0$.  \\
\hspace*{13cm} $\Box$  \\

The kernel of each $\delta$ in the first case is equal to the image of the preceding $\delta$ but this may no longer be true in the restricted case and we set:\\

\noindent
{\bf DEFINITION 3.1.10}: Let $B^s_{q+r}(g_q)\subseteq Z^s_{q+r}(g_q)$ and $H^s_{q+r}(g_q)=Z^s_{q+r}(g_q)/B^s_{q+r}(g_q)$ with $H^1(g_q)=H^1_q(g_q)$ be the coboundary space $im(\delta)$, cocycle space $ker(\delta)$ and cohomology space at ${\wedge}^sT^*\otimes g_{q+r}$ of the restricted $\delta$-sequence which only depend on $g_q$ and may not be vector bundles. The symbol $g_q$ is said to be s-{\it acyclic} if $H^1_{q+r}=...=H^s_{q+r}=0, \forall r\geq 0$, {\it involutive} if it is n-acyclic and {\it finite type} if $g_{q+r}=0$ becomes trivially involutive for r large enough. In particular, if $g_q$ is involutive {\it and} finite type, then $g_q=0$. Finally, $S_qT^*\otimes E$ is involutive for any $ q\geq 0$ if we set $S_0T^*\otimes E=E$. \\

We have (See [17] for the diagram allowing to prove this delicate result first found by Spencer):  \\

\noindent
{\bf PROPOSITION 3.1.11}: If $g_q$ is $2$-acyclic and $g_{q+1}$ is a vector bundle over ${\cal{R}}_q $, then $g_{q+r}$ is a vector bundle over ${\cal{R}}_q, \forall r\geq1$. \\

\noindent
{\bf LEMMA 3.1.12}: If $g_q$ is involutive and $g_{q+1}$ is a vector bundle over ${\cal{R}}_q $, then $ g_q$ is also a vector bundle over 
${\cal{R}}_q$. In this case, changing linearly the local coordinates if necessary, we may look at the maximum number $\beta$ of equations that can be solved with respect to $v^k_{n...n}$ and the intrinsic number $\alpha=m-\beta$ indicates the number of $y$ that can be given arbitrarily.  \\

We notice that ${\cal{R}}_{q+r+1}={\rho}_r ({\cal{R}}_{q+1})$ and ${\cal{R}}_{q+r}=  {\rho}_r({\cal{R}}_q)$ in the following commutative diagram:  \\
\[  \begin{array}{ccc}
             {\cal{R}}_{q+r+1}  &  \stackrel{{\pi}^{q+r+1}_{q+1}}{\longrightarrow }& {\cal{R}}_{q+1}   \\
          \hspace{10mm}   \downarrow   {\pi}^{q+r+1}_{q+r}  &         &  \hspace{7mm}\downarrow {\pi}^{q+1}_q \\
             {\cal{R}}^{(1)}_{q+r}    &  \stackrel{{\pi}^{q+r}_q }{\longrightarrow} & {\cal{R}}^{(1)}_q    \\
               \cap   &   &   \cap   \\
               {\cal{R}}_{q+r} &  \stackrel{{\pi}^{q+r}_q}{\longrightarrow}& {\cal{R}}_q
               \end{array}  \]
but we only have in general ${\cal{R}}^{(1)}_{q+r} \subseteq {\rho}_r({\cal{R}}^{(1)}_q)$. We finally obtain the following crucial Theorem and its Corollary (Compare to [17], p 70-75):  \\

\noindent
{\bf THEOREM 3.1.13}: Let ${\cal{R}}_q\subset J_q({\cal{E}})$ be a system of order $q$ on ${\cal{E}}$ such that ${\cal{R}}_{q+1}$ is a fibered submanifold of $J_{q+1}({\cal{E}})$. If $g_q$ is $2$-acyclic and $g_{q+1}$ is a vector bundle over ${\cal{R}}_q$, then we have ${\cal{R}}^{(1)}_{q+r}={\rho}_r ({\cal{R}}^{(1)}_q)$ for all $r\geq 0$.  \\

\noindent
{\bf DEFINITION 3.1.14}: A system ${\cal{R}}_q\subset J_q({\cal{E}})$ is said to be {\it formally integrable} at the order $q+r$ if ${\pi}^{q+r+s}_{q+r}: {\cal{R}}_{q+r+s}\rightarrow {\cal{R}}_{q+r}$ is an epimorphism of fibered manifolds for all $s\geq 1$, {\it formally integrable} if ${\pi}^{q+r+1}_{q+r}$ is an epimorphism of fibered manifolds $\forall r\geq 0$ and {\it involutive} if it is formally integrable with an involutive symbol $g_q$. We have the following useful test ([17],[43]): \\ 

\noindent
{\bf COROLLARY 3.1.15}: Let ${\cal{R}}_q\subset J_q({\cal{E}})$ be a system of order $q$ on ${\cal{E}}$ such that ${\cal{R}}_{q+1}$ is a fibered submanifold of $J_{q+1}({\cal{E}})$. If $g_q$ is $2$-acyclic (involutive) and if the map ${\pi}^{q+1}_q:{\cal{R}}_{q+1} \rightarrow {\cal{R}}_q$ is an epimorphism of fibered manifolds, then ${\cal{R}}_q$ is formally integrable (involutive).  \\

This is all what is needed in order to study nonlinear systems of ordinary differential (OD) or partial differential (PD) equations, using calligraphic letters like $\cal{E}$ for the nonlinear framework and capital letters like $E=V({\cal{E}})$ for the linear or vertical linearized framework.  \\

\noindent
{\bf 3.2) DIFFERENTIAL ALGEBRA}\\

Let $A, B,\dots$ be commutative unitary rings or even integral domain with fields of quotients $Q(A)=K,Q(B)=L,\dots$ containing a field $k$ as a subring or subfield, for example a polynomial ring $k[x]$ in many indeterminates with coefficients in $k$ and the corresponding field $k(x)$ of rational functions. The ideas that led Erich K\"{a}hler to the next definitions in 1930 are of two kinds ([7]):  \\

\noindent
$\bullet$ The derivative of a polynomial with respect to any one of the indeterminates is a polynomial while the derivative of a rational function is a rational function, a reason sufficient for believing that the concept of derivation could be useful in algebra. \\

\noindent
$\bullet$  The variational and linearization procedures presented in the last section and used for many applications to physics should be extended to differential algebra in order to obtain the algebraic counterpart of  definition 3.1.3 and proposition 3.1.5, replacing $ X$ by a ring $A$ or a field $K$.  \\

\noindent
{\bf DEFINITION 3.2.1}: A {\it derivation} from $A$ to an $A$-module $M$ over $k$ is a map $\delta : A \rightarrow M$ such that $\delta (a+b)=\delta a + \delta b, \delta (ab)=a\delta b + b \delta a$ with $\delta {\mid}_k=0$ and the set of such maps is denoted by $der_k(A,M)$ with $1 \times 1= 1 \Rightarrow \delta 1=0$. When $M=A$, we simply set $ der_k(A,A)=der_k(A)$. \\

\noindent
{\bf PROPOSITION 3.2.2}: Given any $A$-module $M$ and any derivation $\delta \in der_k(A,M)$, there exists a unique $A$-module denoted by ${\Omega}_{A/k} $, called {\it module of K\"{a}hler differentials} of $A$ over $k$, a derivation $d\in der_k(A,{\Omega}_{A/k})$ and a unique morphism $f \in hom_A({\Omega}_{A/k},M)$ such that $\delta = f \circ d $ in the following commutative diagram:  \\
\[  \begin{array}{rcl}
  A & \stackrel{d}{\longrightarrow} & {\Omega}_{A/k}  \\
       &  \delta \searrow \,\, & \,\, \downarrow f   \\
          &  & M
          \end{array}  \]
The element $da\in {\Omega}_{A/k}$ is called the {\it differential} of $a$ and $der_k(A,M)=hom_A({\Omega}_{A/k},M)$.   \\
          
\noindent
{\it Proof}: Let $F$ be the free $A$-module made by the symbols $da, a \in A$ and let $N$ be the submodule of $F$ generated by $d\alpha, d(a+b) - da - db, d(ab) - adb-bda$ for $\alpha \in k, a,b\in A$. We set ${\Omega}_{A/k}=F/N$ and the derivation $d=d_{A/k}: A \rightarrow {\Omega}_{A/k}:a \rightarrow da$ is the {\it universal derivation} allowing to define $f: {\Omega}_{A/k} \rightarrow M$ by $(f\circ d)(a)=f(da)=\delta a$.  \\
Another way is to take into account the limit procedure that is classically used in analysis, namely $(x+h) - x = h \Rightarrow (x+h)^2 - x^2 =2xh + h^2 \Rightarrow dx^2=2xdx$ when $h\mapsto 0$ in order to avoid the square quatity $h^2$ and so on. For this, let us denote by $I$ the kernel of the map $A {\otimes}_k A \rightarrow A: a\otimes b \rightarrow ab$ and {\it define} ${\Omega}_{A/k}=I/I^2$ while setting $da= 1 \otimes a - a \otimes  1, \forall a \in A$. Using the bimodule structure of ${ }_AA_A$ while identifying $a\in A$ with $a\otimes 1 \in A{\otimes}_kA$, it follows that $d$ is indeed a derivation from $A$ to the $A$-module ${\Omega}_{A/k}$ as we have successively:  \\
\[  \begin{array}{rcl}
adb + bda & = & a(1 \otimes b - b \otimes 1) + b(1 \otimes a - a \otimes 1)  \\
   &   =   & (a \otimes 1)(1 \otimes b -b \otimes 1) + (b \otimes 1)(1 \otimes a - a \otimes 1)   \\
    &  =  &  (1 \otimes a)(1 \otimes b - b \otimes 1) + (b \otimes 1)(1 \otimes a - a \otimes 1) - (1 \otimes a - a \otimes 1)(1 \otimes b - b \otimes 1)  \\
      &  =  &  (1 \otimes ab - ab \otimes 1) \hspace{1cm}  mod(I^2)  \\
         &  =  &  d(ab) \hspace{1cm}  mod(I^2)
\end{array}   \]
and thus:   \\
\[  \begin{array}{rcl}
da^2 = 1 \otimes a^2 - a^2 \otimes 1 &  =  &  (1 \otimes a + a \otimes 1 )(1 \otimes a - a \otimes 1)   \\ 
       &   =    &  (a \otimes 1 )(1 \otimes a - a \otimes 1) + (1 \otimes a )(1\otimes a - a \otimes 1 )     \\ 
      &    =    &   2 (a \otimes 1)(1 \otimes a - a \otimes 1)  + (1 \otimes a - a \otimes 1 ) (1 \otimes a - a \otimes 1)  \\
       &  =  &     2 (a \otimes 1 )( 1 \otimes a - a \otimes 1 )  \hspace{1cm} mod(I^2)    \\
          &  =  &  2 a da   \hspace{1cm}   mod (I^2)
\end{array}   \]
\hspace*{13cm}    $\Box$  \\

Among the elementary properties of the K\"{a}hler differentials, we notice that, if $f:A \rightarrow B$ is a $k$-algebra homomorphism and $M$ is a $B$-module, then $M$ becomes a $A$-module under the rule $f(a)=b \Rightarrow am=f(a)m, \forall a \in A, b \in B, m \in M$ and we have the exact sequence of $B$-modules:  \\
\[   0 \rightarrow der_A(B,M) \rightarrow der_k(B,M) \rightarrow der_k(A,M)   \]
because $\delta \in der_k(B,M) \Rightarrow \delta \circ f \in der_k(A,M)$ and thus $der_k(A,M)\simeq hom_B(B {\otimes }_A {\Omega}_{A/k}, M)$.  \\

The proof of the following two propositions is classical and can be found in ([21], p 387-389):  \\

\noindent
{\bf PROPOSITION 3.2.3}: ({\it First fundamental exact sequence}) We have the exact sequence of $B$-modules:   
\[  B {\otimes}_A {\Omega}_{A/k} \rightarrow {\Omega}_{B/k} \rightarrow {\Omega}_{B/A} \rightarrow 0     \]
where the first map is a monomorphism when $f$ is a monomorpism. In particular, if $k \subset K \subset L $ is a chain of field extensions, then one has 
the short exact sequence of vector spaces over $L$:  
\[   0 \rightarrow L{\otimes}_K{\Omega}_{K/k} \rightarrow {\Omega}_{L/k} \rightarrow {\Omega}_{L/K} \rightarrow 0   \]

\noindent
{\bf PROPOSITION 3.2.4}: ({\it Second fundamental exact sequence}) If  we have the short exact sequence $0 \rightarrow \mathfrak{a} \rightarrow A \stackrel{f}{\longrightarrow} B \rightarrow 0$, then we have the exact sequence of $B$-modules:  \\
\[    \mathfrak{a}/{\mathfrak{a}}^2 \rightarrow B {\otimes}_A {\Omega}_{A/k} \rightarrow {\Omega}_{B/k} \rightarrow 0  \]

\noindent
{\bf EXAMPLE 3.2.5}: Let $x,y$ be two indeterminates over the field $k = \mathbb{Q}$ and consider the case $A=k[x,y], B=A/(x^2, y^2)$. Then the ideal $\mathfrak{a}=(x^2,y^2)\subset A$ is not prime because $rad(\mathfrak{a})=\mathfrak{p}=(x,y)\in spec(A)$. The image of $x^3$ is $1 \otimes 3x^2dx$ in $B {\otimes}_A{\Omega}_{A/k}$, that is $3x^2\otimes dx=0$ because $x^2 \in \mathfrak{a}$. A similar comment applies to $y^3$ and it is easy to see that the kernel of the map $\mathfrak{a}/{\mathfrak{a}}^2 \rightarrow B {\otimes}_A{\Omega}_{A/k}$ is of the form $ax^3 + by^3, \forall a,b\in A$.   \\

Finally, if $S$ is a multiplicatively closed subset of $A$, we may use the morphism ${\theta}_S : A \rightarrow S^{-1}A$ in Proposition 2.2.3, we shall study the behaviour of derivations and differentials under localization. As $\delta \in der_k(A,M)$ induces a unique derivation $\delta \in der_k(S^{-1}A, S^{-1}M$ through the known formula 
$\delta (a/s): s\delta a - a \delta s)/s^2$, it follows that the morphism $der_k(S^{-1}A, S^{-1}M) \rightarrow der_k(A, S^{-1}M)$ given by $\delta \rightarrow \delta \circ {\theta}_S$ is an epimorphism. We obtain the short exact sequence:  \\
\[           0 \rightarrow der_A(S^{-1}A,S^{-1}M) \rightarrow der_k(S^{-1}A,S^{-1}M) \rightarrow der_k(A,S^{-1}M) \rightarrow 0  \]
and thus the short exact sequence:  \\
\[     0 \rightarrow S^{-1}A {\otimes}_A{\Omega}_{A/k} \rightarrow {\Omega}_{S^{-1}A/k} \rightarrow {\Omega}_{S^{-1}A/A}\rightarrow 0  \]
Taking into account the previous standard formula, it follows that ${\Omega}_{S^{-1}A/A}=0$ and we obtain:  \\

\noindent
{\bf PROPOSITION 3.2.6}: There is an isomorphism $ S^{-1}A {\otimes}_A{\Omega}_{A/k} \simeq {\Omega}_{S^{-1}A/k}$.  \\

\noindent
{\bf EXAMPLE 3.2.7}:  We now present in an independent manner a few OD or PD cases showing the difficulties met when studying differential ideals and ask the reader to revisit them later on while reading the main Theorems. As only a few results will be proved, the interested reader may look at [18] or [20] for more details and compare to [9] or [11]. \\

\noindent
$\bullet \,\, OD \,\,1$: If $k=\mathbb{Q}$, $y$ is a differential indeterminate and $d_x$ is a formal derivation, we may set $d_xy=y_x, d_xy_x=y_{xx}$ and so on in order to introduce the differential ring $A=k[y,y_x, y_{xx}, ...]=k\{y\}$. We consider the differential ideal $\mathfrak{a}\subset A$ generated by the differential polynomial $P=y^2_x-4y$. We have $d_xP=2y_x(y_{xx}-2)$ and $\mathfrak{a}$ cannot 
be a prime differential ideal, $\dots$ and so on. After no less than $4$ differentiations, we let the reader discover that $y_{xxx}^5 \in \mathfrak{a} \Rightarrow y_{xxx} \in rad(\mathfrak{a})$ and thus $\mathfrak{a}$ is neither prime nor perfect, that is equal to its radical, but $rad(\mathfrak{a})$ is perfect as it is the intersection of the prime differential ideal generated by $y$ with the prime differential ideal generated by $y_x^2-4y$ and $y_{xx}-2$, both containing $y_{xxx}$.  \\

\noindent
$\bullet \,\, OD \,\, 2$: With the same notations, let us consider the differential ideal $\mathfrak{a}\subset A$ generated by the differential polynomial $P=y^2_x-4y^3$. We have $d_xP=2y_x(y_{xx}-6y^2)$ and $\mathfrak{a}$ cannot be  prime differential ideal. Hence, we must have either $y_x=0 \Rightarrow  y=0$ or $y_{xx}-6y^2=0$ and so on. After 3 differentiations  we obtain $(y_{xx}-6y^2)^4 \in \mathfrak{a}  \Rightarrow   y_{xx}-6y^2 \in rad(\mathfrak{a})$ and thus $\mathfrak{a}$ is neither prime nor perfect as before but $rad(\mathfrak{a})$ is the prime differential ideal generated by $y_x^2-4y^3$ and 
$y_{xx}-6y^2$.  \\

\noindent
$\bullet \,\, PD \,\, 1$: If $k=\mathbb{Q}$ as before, $y$ is a differential indeterminate and $(d_1,d_2)$ are two formal derivations, let us consider the differential ideal generated by $P_1=y_{22}-\frac{1}{2}(y_{11})^2$ and $P_2=y_{12}-y_{11}$ in $k\{y\}$. Using crossed derivatives and differentiating twice, we get 
$(y_{111})^3 \in \mathfrak{a} \Rightarrow y_{111} \in rad(\mathfrak{a})$ and thus $\mathfrak{a}$ is again neither prime nor perfect but $rad(\mathfrak{a})$ is a perfect differential ideal and even a prime differential ideal $\mathfrak{p}$ because we obtain easily from the last subsection that the resisual differential ring $k\{y\}/\mathfrak{p}\simeq k[y,y_1,y_2,y_{11}]$ is a differential integral domain. Its quotient field is thus the differential field $K=Q(k\{y\}/\mathfrak{p})\simeq k(y,y_1,y_2,y_{11})$ with the rules: 
\[   d_1y=y_1,d_1y_1=y_{11}, d_1y_{11}=0, d_2y=y_2, d_2y_1=y_{11}, d_2y_{11}=0\] 
as a way to avoid " {\it looking for solutions} ". The formal linearization is the linear system $R_2\subset J_2(E)$ obtained in the last section where it was defined over ${\cal{R}}_2$, but {\it not} over $K$, by the two linear second order PDE:  \\
\[             {\eta}_{22}-y_{11}{\eta}_{11}=0, \hspace{1cm} {\eta}_{12}-{\eta}_{11}=0   \, \,\, \Rightarrow \,\, \,
(y_{11} - 1){\eta}_{111}=0   \] 
changing slightly the notations with $\eta=\delta y$ and keeping the letter $v$ only when looking at the symbols. It is at this point that {\it the problem starts} because ${\cal{R}}_2$ is indeed a fibered manifold with arbitrary parametric jets $(y,y_1,y_2,y_{11})$ but ${\cal{R}}_3={\rho}_1({\cal{R}}_2)$ is no longer a fibered manifold because the dimension of its symbol changes when $y_{11}=1$. We understand therefore that {\it there should be a close link existing between formal integrability and the search for prime differential ideals or differential fields}. The solution of this problem has been provided as early as in 1983 for studying the "Differential Galois Theory " in ([18]). The idea is to add the third order PDE $y_{111}=0$ and thus consider the linearized PDE ${\eta}_{111}=0$ obtaining therefore a third order involutive system well defined over $K$ with symbol $g_3=0$. \\

\noindent
$\bullet \,\, PD \,\, 2$: With the same notations, let us consider the differential ideal generated by the differential polynomials $P_1=y_{22}-\frac{1}{3}(y_{11})^3$ and $P_2=y_{12}-\frac{1}{2}(y_{11})^2$ in $k\{y\}$. We get:  \\
\[ P_1,P_2 \in \mathfrak{a} \Rightarrow d_2P_2 - d_1P_1 + y_{11}d_1P_2=0 \Rightarrow {\cal{R}}_2  \hspace{3mm}involutive \]
with $dim(g_{q})=1, \forall q\geq 1$. As the symbol $g_2$ is involutive, there is an infinite number of parametric jets $(y,y_1,y_2, y_{11}, y_{111}, ...)$ and thus $k\{y\}/\mathfrak{a}\simeq k[y,y_1,y_2, y_{11}, y_{111}, ...]$ is a differential integral domain with $d_2y_2=y_{22}=\frac{1}{3}(y_{11})^3, d_2y_{11}=y_{112}=  y_{11}y_{111}, ...$. It follows that $\mathfrak{a}=\mathfrak{p}$ is a prime differential ideal with $rad(\mathfrak{p})=\mathfrak{p}$. The second order linearized system is:  \\ 
\[  {\eta}_{22} - (y_{11})^2{\eta}_{11}=0, \hspace{1cm}  {\eta}_{12} - y_{11}{\eta}_{11}=0   \]                                                                                                                                                                                                                                                                                                                                                                  
is now well defined over the differential field $K=Q(k\{y\}/\mathfrak{p})$ and is involutive.  \\

\noindent
{\bf DEFINITION 3.2.8}: A {\it differential ring} is a ring $A$ with a finite number of commuting derivations $({\partial}_1, ..., {\partial}_n)$ that can be extended to derivations of the ring of quotients $Q(A)$ as we already saw. We shall suppose from now on that $A$ is even an integral domain and introduce the differential field $K=Q(A)$. For example, if $x^1, ... , x^n$ are indeterminates over $\mathbb{Q}$, then $\mathbb{Q}[x]=\mathbb{Q}[x^1, ...,x^n]$  is a differential ring with quotient differential field $\mathbb{Q}(x)$.  \\

If $K$ is a differential field as above and $(y^1,...,y^m)$ are indeterminates over $K$, we transform the polynomial ring $K\{y\}={lim}_{q\rightarrow \infty}K[y_q]$ into a differential ring by introducing as usual the {\it formal derivations} $d_i={\partial}_i+y^k_{\mu+1_i}\partial/\partial y^k_{\mu}$ and we shall set $K<y>= Q(K\{y\})$.  \\

\noindent
{\bf DEFINITION 3.2.9}: We say that $\mathfrak{a}\subset K\{y\}$ is a {\it differential ideal} if it is stable by the $d_i$, that is if  $d_ia\in\mathfrak{a}, \forall a \in \mathfrak{a}, \forall i=1,...,n$. We shall also introduce the {\it radical} $rad(\mathfrak{a})=\{a\in A\mid \exists r,a^r\in \mathfrak{a}\}\supseteq \mathfrak{a}$ and say that $\mathfrak{a}$ is a {\it perfect} (or {\it radical}) differential ideal if $rad(\mathfrak{a})=\mathfrak{a}$. If $S$ is any subset of $A$, we shall denote by $\{S\}$ the differential ideal generated by $S$ and introduce the (non-differential) ideal ${\rho}_r(S)=\{d_{\nu}a \mid a\in S, 0 \leq  \mid\nu\mid \leq r\}$ in $A$.   \\

\noindent
{\bf LEMMA 3.2.10}: If $\mathfrak{a}\subset A$ is differential ideal, then $rad(\mathfrak{a}) $ is a differential ideal containing 
$\mathfrak{a}$.  \\

\noindent
{\it Proof}: If $d$ is one of the derivations, we have $ a^{r-1}da=\frac{1}{r}da^r \in \{a^r\}$ and thus:  \\
\[ (r-1) a^{r-2}(da)^2 + a^{r-1}d^2a \in \{a^r\}\Rightarrow a^{r-2}(da)^3\in \{a^r\},... \Rightarrow (da)^{2r-1} \in \{a^r\}   \]
\hspace*{13cm}   $\Box$   \\

\noindent
{\bf LEMMA 3.2.11}: If $\mathfrak{a} \subset K\{y\}$, we set ${\mathfrak{a}}_q= \mathfrak{a}\cap K[y_q]$ with ${\mathfrak{a}}_0=\mathfrak{a}\cap K[y]$ and  ${\mathfrak{a}}_{\infty}= \mathfrak{a}$. We have in general ${\rho}_r({\mathfrak{a}}_q) \subseteq {\mathfrak{a}}_{q+r}$ and the problem will be to know when we may have equality.  \\ 

We shall say that a differential extension $L=Q(K\{y\}/\mathfrak{p})$ is a {\it finitely generated} differential extension of $K$ and we may define the {\it evaluation epimorphism} $K\{y\} \rightarrow K\{\eta \}\subset L$ with kernel $\mathfrak{p}$ by calling $\eta$ or $\bar{y}$ the residue of $y$ modulo $\mathfrak{p}$. If we study such a differential extension $L/K$, by analogy with Section 2, we shall say that $R_q$ or $g_q$ is a vector bundle over ${\cal{R}}_q$ if one can find a certain number of maximum rank determinant $D_{\alpha}$ that cannot be all zero at a generic solution of ${\mathfrak{p}}_q$ defined by differential polynomials $P_{\tau}$, that is to say, according to the Hilbert Theorem of Zeros, we may find polynomials $A_{\alpha}, B_{\tau}\in K\{y_q\}$ such that $  {\sum}_{\alpha} A_{\alpha}D_{\alpha} + {\sum}_{\tau} B_{\tau}P_{\tau}=1 $. The following Lemma will be used in the next important Theorem:  \\

\noindent
{\bf LEMMA 3.2.12}: If $\mathfrak{p}$ is a prime differential ideal of $K\{y\}$, then, for $q$ sufficiently large, there is a polynomial $D\in K[y_q]$ such that $D\notin {\mathfrak{p}}_q$ and :   \\
\[          D{\mathfrak{p}}_{q+r} \subset rad ({\rho}_r({\mathfrak{p}}_q)) \subset {\mathfrak{p}}_{q+r}, \hspace {1cm}  \forall r\geq 0  \] 

\noindent
{\bf THEOREM  3.2.13}: ({\it Primality test}) Let ${\mathfrak{p}}_q\subset K[y_q]$ and ${\mathfrak{p}}_{q+1}\subset K[y_{q+1}]$ be prime ideals  such that ${\mathfrak{p}}_{q+1}={\rho}_1({\mathfrak{p}}_q)$ and ${\mathfrak{p}}_{q+1}\cap K[y_q]={\mathfrak{p}}_q$. If the symbol $g_q$ of the algebraic variety ${\cal{R}}_q$ defined by ${\mathfrak{p}}_q$ is $2$-acyclic and if its first prolongation $g_{q+1}$ is a vector bundle over ${\cal{R}}_q$, then $\mathfrak{p}={\rho}_{\infty}({\mathfrak{p}}_q)$ is a prime differential ideal with $\mathfrak{p} \cap K[y_{q+r}]={\rho}_r({\mathfrak{p}}_q), \forall r\geq 0 $.  \\

\noindent
{\bf COROLLARY  3.2.14}: Every perfect differential ideal of $\{y\}$ can be expressed in a unique way as the non-redundant intersection of a finite number of prime differential ideals.  \\

\noindent
{\bf COROLLARY 3.2.15}: ({\it Differential basis}) If $\mathfrak{r}$ is a perfect differential ideal of $K\{y\}$, then we have $\mathfrak{r}=rad({\rho}_{\infty} ({\mathfrak{r}}_q))$ for $q$ sufficiently large.  \\

\noindent
{\bf EXAMPLE 3.2.16}: As $K\{y\}$ is a polynomial ring with an infinite number of variables it is not noetherian and an ideal may not have a finite basis. With $K=\mathbb{Q}, n=1$ and $d=d_x$, then $\mathfrak{a}=\{yy_x,y_xy_{xx},y_{xx}y_{xxx}, ... \} \Rightarrow (y_x)^2+yy_{xx}\in \mathfrak{a} \Rightarrow rad(\mathfrak{a})=\{y_x\}$ is a prime differential ideal.  \\
 
\noindent
{\bf PROPOSITION  3.2.17}: If $\zeta$ is differentially algebraic over $K<\eta>$ and $\eta$ is differentially algebraic over $K$, then $\zeta$ is differentially algebraic over $K$. Setting $\xi=\zeta - \eta$, it follows that, if $L/K$ is a differential extension and $\xi,\eta \in L$ are both differentially algebraic over $K$, then $\xi + \eta$, $\xi\eta$ and $d_i\xi$ are differentially algebraic over $K$.  \\

If $L=Q(K\{y\}/\mathfrak{p})$, $M=Q(K\{z\}/\mathfrak{q})$ and $N=Q(K\{y,z\}/\mathfrak{r})$ are such that $\mathfrak{p}=\mathfrak{r}\cap K\{y\}$ and $\mathfrak{q}=\mathfrak{r}\cap K\{z\}$, we have the two towers $K\subset L\subset N$ and $K\subset M\subset N$ of differential extensions and we may therefore define the new tower $K \subseteq L\cap M \subseteq <L,M> \subseteq N$. However, if only $L/K$ and $M/K$ are known and we look for such an $N$ containing both $L$ and $M$, we may use the universal property of tensor products an deduce the existence of a differential morphism $L{\otimes}_KM\rightarrow N$ by setting $d(a\otimes b)=(d_La) \otimes b+a \otimes (d_Mb)$ whenever $d_L\mid K=d_M\mid K=\partial$. Looking for an abstract {\it composite} differential field amounts therefore to look for a prime differential ideal in $L{\otimes}_K M$ which is a direct sum of integral domains (See [18] for more details).  \\

\noindent
{\bf DEFINITION  3.2.18}: A differential extension $L$ of a differential field $K$ is said to be {\it differentially algebraic} over $K$ if every element of $L$ is differentially algebraic over $K$. The set of such elements is an intermediate differential field $K' \subseteq L$, called the {\it differential algebraic closure} of $K$ in $L$. If $L/K$ is a differential extension, one can always find a maximal subset $S$ of elements of $L$ that are differentially transcendental over $K$ and such that $L$ is differentially algebraic over $K<S>$. Such a set is called a {\it differential transcendence basis} and the number of elements of $S$ is called the {\it differential transcendence degree} of $L/K$.  \\  

\noindent
{\bf THEOREM 3.2.19}: If $L/K$ is a finitely generated differential extension, then any intermediate differential field $K'$ between $K$ and $L$ is also finitely generated over $K$.  \\

\noindent
{\bf THEOREM  3.2.20}: The number of elements in a differential basis of $L/K$ does not depent on the generators of $L/K$ and his value 
is $difftrd(L/K)=\alpha$. Moreover, if $K\subset L \subset M$ are differential fields, then $difftrd(M/K)=difftrd(M/L) + difftrd(L/K)$.  \\

Comparing the differential geometric approach to nonlinear algebraic systems with the differential algebraic approach just presented while setting $\delta y_q={\eta}_q$, we obtain: \\

\noindent
{\bf COROLLARY 3.2.21}: When $L/K$ is a finitely generated differential extension, then ${\Omega}_{L/K}$ is a differential module over the differential ring 
$L {\otimes}_KK[d]=L[d]$ of differential operators with coefficients in $L$. The linearized "{\it system} "  $R=hom_L({\Omega}_{L/K}, L)$ is thus a (left) differential module for the Spencer operator like in the linear framework.  \\

It is not evident to grasp these results in order to apply them to control theory or mathematical physics for two reasons. The first is that the formal theory of nonlinear systems has not been accepted by differential geometers because of the homological background based on the so-called "{\it vertical machinery}" and the systematic use of the Spencer $\delta$-cohomology. The recent study of the Schwarzschild and Kerr metrics (Compare [35] to [1]) is providing a good example of such a poor situation. The second is the fact that, when $K$ is a true differential field and $M$ is differential module defined over the noncommutative ring $K[d_1,\dots,d_n]=K[d]$ of differential operators with coefficients in $K$, then the "{\it system}" $R=hom_K(M,K)$ is still not used today because its differential structure highly depends on the Spencer operator which has never been introduced in physics. As a good example, we may quote the fact that the {\it Cosserat couple stress equations} are just described by the formal adjoint of the linear Spencer operator ([19]). \\

\noindent
{\bf 3.3) NONLINEAR CONTROL THEORY}  \\

As we have already explained in ([18]), the generalized " {\it B\"{a}cklund problem} " is nothing else than the study of nonlinear {\it differential correspondences} in the theory of differential elimination. We shall provide, below and successively, a differential geometric definition followed by a differential algebraic definition and all the problem will be to establish a link between them.   \\

When $X$ is a manifold of dimension $n$, let us consider two fibered manifolds over $X$, namely 
${\cal{E}}$ with lcal coordinates $(x,y)$ and ${\cal{F}}$ with local coordinates $(x,z)$. The fibered roduct ${\cal{E}} {\times}_X {\cal{F}}$ is a fibered manifold over $X$ with ocal coordinates $(x,y,z)$ and we have the canonical identification:   \\
\[   J_q({\cal{E}} {\times}_X {\cal{F}}) = J_q({\cal{E}}) {\times}_XJ_q({\cal{F}})    \]
with local coordinates $(x,y_q,z_q)$.   \\
For most applications, we shall suppose that ${\cal{E}}=X \times Y$ and ${\cal{F}}=X \times Z$.  \\

\noindent
{\bf DEFINITION 3.3.1}: Let ${\cal{R}}_q \subset J_q({\cal{E}} {\times}_X {\cal{F}})$ be a nonlinear system of order $q$ on ${\cal{E}}{\times}_X {\cal{F}}$ called a {\it differential correspondence} between $(y,z)$. When $r \rightarrow \infty$, we may consider the {\it resolvent systems} ${\cal{P}}_{q+r} \subset J_{q+r}({\cal{E}})$ for $y$ and ${\cal{Q}}_{q+r} \subset J_{q+r}({\cal{F}})$ for $z$, induced by the canonical projections of ${\cal{E}} {\times}_X {\cal{F}}$ onto ${\cal{E}}$ and ${\cal{F}}$ respectively.    \\

Roughly, finding ${\cal{P}}$ amounts to eliminate $z$ while finding ${\cal{Q}}$ amounts to eliminate $y$ and we shall only consider te first problem as the second will be similar.  \\

\noindent
$\bullet$  {\it In the linear case}, pushing $y$ on the left and $z$ on the right, we are left with the search of the CC for $y$ or the CC for $z$ that may be quite difficult. One of the best examples has been provided by M. Janet with the second order system (See [20] or [22]  for details):  \\
\[         y_{33} - x^2 y_{11}=z^1, \hspace{2cm}   y_{22}=z^2     \]
where $y=f(x)$ can be given arbitrarily for getting $z=g(x)$ while $z=g(x)$ must satisfy one CC of order $3$ and one CC of order $6$.  \\ 

\noindent
$\bullet$  {\it In the nonlinear case}, we have ([18],[21]):    \\

\noindent
{\bf THEOREM OF THE RESOLVENT SYSTEMS 3.3.2}: In general, one may find two integers $r,s \geq 0$ such that $ {\cal{R}}^{(s)}_{q+r}$ is formally integrable (involutive) with formally integrable (involutive) projections ${\cal{P}}^{(s)}_{q+r} \subset J_{q+r}({\cal{E}})$ and ${\cal{Q}}^{(s)}_{q+r} \subset J_{q+r}({\cal{F}})$. Moreover, $r$ and $s$ can be (tentatively) found by a finite algorithm preserving the symmetry existing between $\cal{E}$ and $\cal{F}$.  \\ 

\noindent
{\it Proof}: First of all, we know that, in general, one can find the two integers $r,s \geq 0$ in such a way that ${\cal{R}}^{(s)}_{q+r}$ is formally integrable (involutive). Hence, using the commutative and exact diagram:  \\

\[   \begin{array}{cccl}
{\cal{R}}_{q+r+s} & \longrightarrow & {\cal{P}}_{q+r+s} & \rightarrow 0  \\
   \downarrow &   &  \downarrow  &     \\
   {\cal{R}}^{(s)}_{q+r} &  \longrightarrow & {\cal{P}}^{(s)}_{q+r} & \rightarrow 0  \\
       \cap   &       &   \cap    &       \\
    {\cal{R}}_{q+r} &  \longrightarrow  &  {\cal{P}}_{q+r}  &  \rightarrow 0
\end{array}   \]
we may suppose, without any loss of generality, that ${\cal{R}}_q$ is formally integrable (involutive). \\
Now, chasing in the commutative diagram:   \\
\[    \begin{array}{ccccccl}
   &  &  J_s({\cal{R}}_{q+r})  &  & \longrightarrow &    & J_s(J_{q+r}({\cal{E}}{\times}_X {\cal{F}})  \\
  & \nearrow   &    |     &        &    &  \nearrow & \hspace{8mm} |   \\
 {\cal{R}}_{q+r+s}  &   & \longrightarrow &  &  J_{q+r+s}({\cal{E}}{\times}_X {\cal{F}})& & \hspace{8mm} |  \\
   |    &    & \downarrow    &      &   |     &  & \hspace{8mm}   \downarrow     \\    
 | &  &  J_s({\cal{P}}_{q+r})   &   &\longrightarrow    &    & J_s(J_{q+r}({\cal{E}}))  \\
 \downarrow         & \nearrow   &    &    &  \downarrow  &  \nearrow  &     \\ 
{\cal{P}}_{q+r+s}     &               & \longrightarrow  & &  J_{q+r+s}({\cal{E}})  &   &
 \end{array}     \]
we obtain therefore ${\cal{P}}_{q+r+s} \subseteq {\rho}_s ({\cal{P}}_{q+r})= J_s({\cal{P}}_{q+r}) \cap J_{q+r+s}({\cal{E}}) \subset J_s(J_{q+r}({\cal{E}})), \forall r,s \geq 0$.    \\
Then, chasing in the commutative diagram:   \\
\[      \begin{array}{cccl}
{\cal{R}}_{q+r+s}&  \longrightarrow  &  {\cal{P}}_{q+r+s} & \rightarrow 0  \\
    \downarrow  &   &  \downarrow  &    \\
    {\cal{R}}_{q+r} &  \longrightarrow & {\cal{P}}_{q+r}  &  \rightarrow 0  \\
    \downarrow &    &   &   \\
    0  &   &   &   
\end{array}   \]
we notice that ${\pi}^{q+r+s}_{q+r}: {\cal{P}}_{q+r+s} \longrightarrow {\cal{P}}_{q+r}$ is an epimorphism $\forall r,s \geq 0$.  \\
Finally, chasing in the commutative and exact diagram:  \\
\[    \begin{array}{rcccl}
0 \rightarrow & {\cal{P}}_{q+r+s} & \longrightarrow & {\rho}_s({\cal{P}}_{q+r}) & \rightarrow 0  \\
       &     \downarrow  &    &   \downarrow  &    \\
 0 \rightarrow & {\cal{P}}_{q+r} &   =  & {\cal{P}}_{q+r}  &   \rightarrow 0  \\
    &  \downarrow &  &  &  \\
    &  0  &   &   &  
  \end{array}   \]
we deduce that {\it each} ${\cal{P}}_{q+r}$ is formaly integrable at {\it each}  $q+r, \forall r \geq 0$, though not always formally integrable as we shall see on examples.  \\
Looking at the symbol $h$ of ${\cal{P}}$, we have $h_{q+r+s} \subseteq {\rho}_s(h_{q+r})$ over ${\cal{P}}_{q+r+s}$. According to standard Noetherian arguments, such a situation is stabilizing for $r$ and $s$ large enough but such an approach is not constructive in general.  \\
 
For this reason, we shall prefer to use a different approach which is closer to the one met in the case of linear differential correspondences. For this, if $z=g(x)$ is an arbitrary section of ${\cal{F}}$, we shall consider the new system for $y$ defined by ${\cal{A}}_q=j_q(f)^{-1}({\cal{R}}_q)$ over $K<g>$. Such a system, which is in general neither involutive nor even formally integrable as we shall see on examples, may also be not even compatible as it may not provide a fibered manifold but this way may give informations on the order of the OD or PD equations that should be satisfied by $z$. A similar procedure could be used by setting $y=f(x)$ and introducing 
${\cal{B}}_q=j_q(f)^{-1}({\cal{R}}_q)$ in order to obtain a system for $z$ over $K<f>$.  \\
\hspace*{13cm}    $\Box$     \\

Let us now turn to the differential algebraic counterpart.   \\  

\noindent
{\bf DEFINITION 3.3.3}: If $K$ is a differential field and we have a {\it differential algebraic correspondence} defined by a prime differential ideal $\mathfrak{r} \subset K \{ y,z \}$, we may define the resolvent system for $y$ by the {\it resolvent differential ideal} $\mathfrak{p}=\mathfrak{r} \cap K\{ y \}$ and the resolvent system for $z$ by the {\it resolvent differential ideal} $\mathfrak{q}=\mathfrak{r} \cap K \{ z \}$.  \\

\noindent
{\bf LEMMA 3.3.4}: The resolvent ideal for $y$ is the prime differential {\it resolvent ideal} 
$\mathfrak{p}= \mathfrak{r} \cap K\{ y \}$ for which one can find a differential basis. Similarly, the prime differential resolvent ideal for $z$ is $\mathfrak{q} = \mathfrak{r} \cap K\{ z \} $.  \\

\noindent
{\it Proof}: We have the commutative and exact diagram:   \\
\[   \begin{array}{rcccccl}
      &   0   & &  0   &  &   0 &   \\
        &  \downarrow &     &  \downarrow &   &  \downarrow  &   \\
   0  \rightarrow &  \mathfrak{p}  &   \longrightarrow  &   K \{ y \}  &  \longrightarrow  &  A   & \rightarrow 0  \\
   &  \downarrow &     &  \downarrow &   &  \downarrow &  \\
0  \rightarrow &  \mathfrak{r}  &  \longrightarrow &  K\{  y,z \} & \longrightarrow & B   & \rightarrow 0  
\end{array}   \]
First of all, $B$ is an integral domain because $\mathfrak{r}$ is a prime differential ideal. It follows from a chase that the induced morphism $A \rightarrow B$ is a monomorphism and $A\simeq im(A) \subset B$ is thus also an integral domain, a result showing that $\mathfrak{p}$ is a prime differential ideal. It is essential to notice that {\it projections of ideals cannot be used in the nonlinear framework}. Hence, the idea is to reduce the study of differential algebraic correspondences to the study of purely algebraic correspondences.  \\
\hspace*{13cm}    $\Box$     \\

We end this last section with a few basic motivating examples showing the importance of the non-commutative localization of integral domains for explicit computations and applications. We hope therefore that these examples could be used as test examples for future applications of computer algebra (Compare to [16]).  \\

\noindent
{\bf EXAMPLE 3.3.5}: With $m=1, n=2, q=2, K=\mathbb{Q}(x^1,x^2)$ while using local coordinates $(x^1, x^2, y)$ for the fibered manifold ${\cal{E}}$ let us consider anew the nice example presented by J. Johnson in ([8]), namely the nonlinear system ${\cal{R}}_2 \in J_2({\cal{E}})$ defined by the two algebraic PD equations :      \\
\[   P_1 \equiv y_{22} - x^2y = 0, \,\,\,\,  P_2 \equiv y_{12} - (y)^2 = 0  \]
We let the reader prove successively as an exercise that:  \\
${\cal{R}}^{(1)}_2$ is adding $ 2 y y_2 - x^2 y_1=0 $.  \\
${\cal{R}}^{(2)}_2$ is adding $ 2 (y_2)^2 - y_1 +x^2 (y)^2 = 0 $.  \\
${\cal{R}}^{(3)}_2$ is adding $6x^2yy_2=0$ and thus $y_1=0$.  \\
${\cal{R}}^{(4)}_2$ is adding $ (y)^2=0 $.  \\
Accordingly, the prime ideal ${\mathfrak{p}}_2 \in K[y,y_1,y_2,y_{11}, y_{12},y_{22}]$ generated by the two given differential polynomias $(P_1, P_2)$ is such that $y\in rad({\rho}_4({\mathfrak{p}}_2)) \Rightarrow rad({\rho}_{\infty}({\mathfrak{p}}_2))=(y)$, a result not evident at first sight and leading to the trivial differential extension 
$L=K$. The linearization procedure is even less evident. Indeed, starting with the linearized second order system:  
\[   \delta P_1 \equiv  {\eta}_{22} - x^2 \eta=0, \,\,\, \delta P_2 \equiv {\eta}_{12}- 2 y \eta=0\]
we let the reader prove that we successively get:  \\
$R^{(1)}_2$ is adding $ 2y {\eta}_2 - x^2 {\eta}_1 + 2y_2 \eta=0$.  \\
$R^{(2)}_2$ is adding $  4y_2 {\eta}_2 - {\eta}_1 + 2x^2 y \eta=0$.  \\
$R^{(3)}_2$ is adding $ {\eta}_1=0$.  \\
$R^{(4)}_2$ is adding $y\eta=0$ but one cannot conclude.  \\
Such an example is proving that, {\it in general}, one must start from a formally integrable or even involutive system in order to be able to define the module of K\"{a}hler differentials for the differential extension $L/K$.  \\

\noindent
{\bf EXAMPLE 3.3.6}:  ({\it Burgers}) With $n=2, m=2$, local coordinates $(x^1,x^2, y,z)$ and differential field $K=\mathbb{Q}$, let us consider the algebraic first order involutive system ${\cal{R}}_1$ defined by two differential algebraic PD equations:   \\
\[   \left\{    \begin{array}{l}
z_2 - y=0  \\
y_2 - z_1 + (y)^2=0 
\end{array}  \right. \fbox { $  \begin{array}{cc}
1& 2 \\
1 & 2 
\end{array}  $ }   \]
These two differential polynomials generate a prime differential ideal $\mathfrak{r} \subset K\{ y,z\}$ and provide thus a differential extension $N/K$. Indeed, $K[y,y_1,y_2;z,z_1,z_2]/{\mathfrak{r}}_1\simeq K[y,y_1;z,z_1]$ is an integral domain and ${\mathfrak{r}}_1$ is a prime ideal. Then, using one prolongation, we may introduce the following second order system ${\cal{R}}_2 = {\rho}_1 ({\cal{R}}_1)$:  \\
\[  \left\{  \begin{array}{l}
z_{22} - y_2 =0  \\
y_{22}+ 2 y y_2 - y_1= 0   \\
z_{12} - y_1 = 0  \\
y_{12} -z_{11} + 2 y y_1 = 0   \\
z_2 - y= 0   \\
y_2 - z_1 + (y)^2  = 0
\end{array} \right.   \fbox {$  \begin{array}{cc}
1& 2  \\
1 & 2  \\
1 & \bullet  \\
1 &  \bullet \\   
\bullet & \bullet  \\
\bullet & \bullet
\end{array}  $ } \]
and use the Janet tabular to prove that it is a nonlinear involutive system. It follows that $K[y, ...,y_{22};z,...,z_{22}]/ {\mathfrak{r}}_2 \simeq K[y,y_1,y_{11}; z, z_1, z_{11}] $ is an integral domain and ${\mathfrak{r}}_2$ is a prime ideal. Thanks to Theorem $3.2.13$, we obtain thus finally $K\{y,z\}/ \mathfrak{r} \simeq K[y,y_1,y_{11}, ...; z, z_1, z_{11}, ...] $ which is also an integral domain. It follows hat $\mathfrak{p}=\mathfrak{r}\cap K\{y\}$ and $\mathfrak{q} = \mathfrak{r} \cap K\{z\}$ are prime differetial ideals.  \\
Taking now {\it any} section $y=f(x) $, we obtain the system ${\cal{B}}_1={j_1(f)}^{-1}({\cal{R}}_1)$ for 
$z$:   \\
\[  \left\{  \begin{array}{l}
z_2 - f(x)= 0   \\
z_1 - {\partial}_2f(x) - (f(x))^2 = 0
\end{array} \right.  \]
and its first prolongation ${\cal{B}}_2={j_2(f)}^{-1}({\cal{R}}_2)$ for $z$:  \\
\[  \left\{  \begin{array}{l}
z_{22} - {\partial}_2 f(x) = 0  \\
z_{12} - {\partial}_1f(x) = 0  \\
z_{11}-{\partial}_{12}f(x) - 2 f(x) {\partial}_1f(x) = 0   \\
z_2 - f(x)= 0   \\
z_1 - {\partial}_2f(x) - (f(x))^2 = 0   \\
{\partial}_{22}f(x)+ 2 f(x) {\partial}_2 f(x)- {\partial}_1f(x) = 0 
\end{array} \right.  \]
First of all, this is a fibered manifold if and only if $f$ is solution of the second order system ${\cal{P}}_2$ defined by the single second order PD equation:    \\
\[   y_{22} + 2 yy_2 - y_1 = 0  \]
which is the resolvent system for $y$ generating the prime differential ideal $\mathfrak{p} \subset K\{y\}$ allowing to define a differential extension $L=Q(K\{y\}/\mathfrak{p})$ of $K$ and we have $\mathfrak{p} = \mathfrak{r} \cap K\{y\}\Rightarrow L \subset N$. \\
We are thus left with the first order (nonlinear) system for $z$:  \\
\[  \left\{  \begin{array}{l}
z_2 - f(x)= 0   \\
z_1 - {\partial}_2f(x) - (f(x))^2 = 0
\end{array} \right.  \fbox { $ \begin{array}{cc}
1 & 2 \\
1& 2 
\end{array}   $  }   \]
which is easily seen to be involutive for any solution $y=f(x) $ of ${\cal{P}}_2$. \\
Taking finally any section $z=g(x)$, we obtain the system ${\cal{A}}_1=j_1(g)^{-1}({\cal{R}}_1)$:  \\
\[    \left \{  \begin{array}{l}
 y_2 + (y)^2 -{\partial}_1g(x)=0\\
 y - {\partial}_2g(x)=0
 \end{array} \right.   \]
and the projecion ${\cal{A}}^{(1)}_1$ of its first prolongation ${\cal{A}}_2={j_2(g)}^{-1}({\cal{R}}_2)$:   \\
\[  \left\{  \begin{array}{l}
y_2  + (y)^2 - {\partial}_1 g = 0  \\
y_1 - {\partial}_{12}g = 0  \\
y -{\partial}_2 g= 0   \\
{\partial}_{22}g(x) + ({\partial}_2g(x))^2 - {\partial}_1g(x)=0
\end{array} \right.  \] 
is compatible if and only if $g$ is solution of the second order system ${\cal{Q}}_2$ defined by the single second order PD equation obtained after substitution of $y={\partial}_2g(x)$:  \\
\[          z_{22} + (z_2)^2  - z_1= 0   \]
which is the resolvent system for $z$ generating the prime differential ideal $\mathfrak{q} = \mathfrak{r} \cap K\{z\}$ allowing to define a differential extension $M=Q(K\{z\}/\mathfrak{q})$ of $K$ and we have $\mathfrak{q}= \mathfrak{r}\cap K\{z\}\Rightarrow M \subset N$. \\
We are thus left with the {\it only} zero order (linear) equation for $y$, namely:   \\
\[     y - {\partial}_2 g(x)= 0  \]
for any solution $z=g(x)$ of ${\cal{Q}}_2$. The differential correspondence that must be used is thus ${\cal{R}}_2$.   \\
Both $L, M$ and $ N$ are differential algebraic extensions of $K$ of zero differential transcendence degree.   \\

\noindent
{\bf EXAMPLE 3.3.7}: ({\it Korteweg-de Vries}) With the same notations, we let the reader provide the details of the following similar example with the second order nonliear differential correspondence ${\cal{R}}_2$:    \\
\[   z_2 + (z)^2 + 2 y = 0, \hspace{5mm}   y_{22 }+ 2 y (z)^2 + 4(y)^2 - 2 z y_2 - \frac{1}{2} z_1 = 0  \]                                                                                                    
by exhibiting the nonlinear formally integrable involutive system ${\cal{R}}^{(1)}_2$.      \\
\[  \left\{  \begin{array}{l}
z_{22} + 2 z z_2 +2y_2=0   \\
y_{22} + 2 y(z)^2 + 4 (y)^2 - 2 z y_2 - \frac{1}{2}z_1 = 0  \\
z_{12} + 2 z z_1+ + 2 y_1 =0 \\
z_2 + (z)^2 + 2 y =0
\end{array}  \right.  \fbox { $ \begin{array}{cc}
1  &  2  \\
1 & 2  \\
1 &  \bullet  \\
\bullet & \bullet
\end{array} $ }  \]
for $(y,z)$ such that ${\cal{R}}^{(1)}_3={\rho}_1({\cal{R}}^{(1)}_2)$ according to Theorem $3.1.13$, with characters ${\alpha}^1_2=3, {\alpha}^2_2=0$ (Compare to [18]).  \\
Taking now {\it any} section $y=f(x)$, we obtain the system ${\cal{B}}_2 = j_2(f)^{-1}({\cal{R}}^{(1)}_2)$ for $z$ which is the first prolongation of the first order (nonlinear) system ${\cal{B}}_1$ defined by: \
\[ \left \{  \begin{array}{l}
z_2 + (z)^2 + 2 f(x) =0  \\
z_1 + 4 {\partial}_2f(x)z - 4 f(x)(z)^2 - 8 (f(x))^2 - 2{\partial}_{22}f(x)=0 
\end{array} \right.  \fbox { $ \begin{array}{cc}
1 & 2 \\
1 & \bullet 
\end{array} $ }  \]
Using crossed derivatives and tedious but elementary substitutions, this system is involutive if and only if $y=f(x)$ is a solution of the third order involutive resolvent system ${\cal{P}}^{(1)}_3$ for $y$:  \\
\[     y_{222} + y_1 + 12 y y_{2} = 0   \]
Similarly, taking now {\it any} section $z=g(x)$, we obtain the system ${\cal{A}}_2= j_2(g)^{-1}({\cal{R}}^{(1)}_2)$ defined by:  \\
\[  \left \{ \begin{array}{l}
2 y_{22} + 4(g(x))^2 + 8 (y)^2 -4 g(x) y_2 - {\partial}_1g(x)=0    \\
2 y_2 + {\partial}_{22}g(x)+ 2 g(x){\partial}_2g(x)=0  \\
2 y_1 + {\partial}_{12}g(x) + 2 g(x) {\partial}_1g(x)=0  \\
2y + {\partial}_2g(x) + (g(x))^2 = 0
\end{array}  \right.  \]
Differentiating the second equation with respect to $x^2$ and substracting the first while using the other equations, we discover that this system is compatible if and only if $z=g(x)$ is a solution of the third order involutive resolvent system ${\cal{Q}}^{(1)}_3$ for $z$:    \\
\[     z_{222} + z_1 - 6 (z)^2 z_2 =0   \]
It follows that we are left with a single zero order equation for $y$, namely:  \\  
\[      2y + {\partial}_2g(x) + (g(x))^2 = 0  \]
for any solution $z=g(x)$ of ${\cal{Q}}_2$. The differential correspondence that must be used is thus ${\cal{R}}^{(1)}_3$.   \\

\noindent
{\bf EXAMPLE 3.3.8}: With $n=1,m=2, K=\mathbb{Q}$, let us consider the single input/single output (SISO) nonlinear control system $P \equiv y^1y^2_x - y^1_x - a=0$ with a constant parameter $a\in K$. The differential ideal $\mathfrak{p} \subset K\{y^1,y^2\}=K\{y\}$ generated by $P$ is prime because $K \{y\}/\mathfrak{p}=K[y^1; y^2,y^2_x, ... ]$ is an integral domain and we set as usual 
$L=Q(K\{y\}/ \mathfrak{p})$. The corresponding linearized system is  
$  y^1 {\eta}^2_x+ y^2_x{\eta}^1 - {\eta}^1_x=0$. Multiplying by a test function $\lambda$ and integrating by parts, the adjoint operator is:   \\
\[      \left \{   \begin{array}{rcl}
{\eta}^1 & \rightarrow &  d\lambda + y^2_x \lambda = {\mu}^1  \\
{\eta}^2 &  \rightarrow & - y^1 d\lambda - y^1_x \lambda = {\mu}^2
\end{array}  \right.  \]
Multiplying the first OD equation by $y^1$, the second by $1$ and adding them, we get 
$a \lambda = y^1 {\mu}^1 + {\mu}^2$. As $L[d]$ is a principal ideal domain, it follows that 
$M={\Omega}_{L/K}$ is a torsion-free and thus free differential module over $L[d]$ if and only if this operator is injective ([21],[22]), that is to say if and only if $a\neq 0$. \\
If $a=0$, then $t(M)$ is generated by $\omega =  y^1{\eta}^2 - {\eta}^1= y^1 \delta y^2 - \delta y^1$ satisfying:  \\
\[    d \omega = y^1 \delta y^2_x + y^1_x \delta y^2 - \delta y^1_x= \delta (y^1y^2_x) + y^2_x \omega - \delta y^1_x = y^2_x \omega \,\,\,  \Rightarrow \,\,\,  d\omega - y^2_x \omega = 0  \]
As $\delta \omega= \delta y^1 \wedge \delta y^2 $, we obtain $\omega \wedge \delta\omega =0$ and one can thus use the analytic Frob\'enius theorem with integrating factor $y^1$ in order to get 
$\omega = y^1\delta (y^2 -log (y^1))$. \\
If $a\neq 0$, say $ a=1$, we have $\lambda= y^1 {\mu}^1 + {\mu}^2$ and obtain the only CC:  \\
\[        y^1 {\mu}^1_x +2 \, y^1_x {\mu}^1 + {\mu}^2_x + y^2_x{\mu}^2 =0  \]   
Multiplying by a test function $\xi$ and integrating by parts, we obtain the parametrization:  \\
\[   - y^1 d{\xi}+ y^1_x{\xi}= {\eta}^1,   \hspace{2cm}    - d {\xi} + y^2_x {\xi}= {\eta}^2  \]
which is injective with potential $\xi= y^1 {\eta}^2 - {\eta}^1$.  \\

\noindent
{\bf EXAMPLE 3.3.9}: With $n=1, m=3, K=\mathbb{Q}$, let us consider the first order nonlinear system ([21]):  \\
\[         P \equiv  2 y^3_x + (y^2_x)^2 - (y^1_x)^2=0   \]
The differenial ideal $\mathfrak{p} $ is prime because we have $K\{y\}/\mathfrak{p}= K[y^1, y^1_x , ...;y^2, y^2_x, ... ;y^3]$ is an integral domain and we define as usual the differential extension $L=Q(L\{y\}/\mathfrak{p})$.      \\
Setting $\delta y = \eta$ and dividing by $2$, the linearized system becomes:  \\
\[           {\eta}^3_x + y^2_x {\eta}^2_x - y^1_x{\eta}^1_x = 0    \]
Multiplying as usual by the Lagrange multiplier $\lambda$ and integrating by parts, we get the adjoint operator with $d=d_x$:  \\
\[      \left \{   \begin{array}{rcl}
{\eta}^1& \rightarrow & y^1_x d\lambda + y^1_{xx} \lambda = {\mu}^1  \\
{\eta}^2&  \rightarrow & - y^2_x d\lambda - y^2_{xx} \lambda = {\mu}^2   \\ 
{\eta}^3&  \rightarrow &  - d\lambda = {\mu}^3
\end{array}  \right.  \]
which is injective with the two CC:    \\
\[     \frac{  {\mu}^1 + y^1_x {\mu}^3}{y^1_{xx}} + \frac{ {\mu}^2 - y^2_x {\mu}^3}{ y^2_{xx}} =0 ,  \hspace{2cm}
d(\frac{{\mu}^1 + y^1_x {\mu}^3}{y^1_{xx}}) + {\mu}^3= 0  \]
It follows that $M={\Omega}_{L/K}$ is a torsion-free differential module over $L[d]$ which is thus also free because it is known that any module over a principal ideal ring which is torsion-free is also free ([39]). Its adjoint operator provides therefore the first order parametrization:  \\
\[   \left \{   \begin{array}{rcl}
{\mu}^1 & \rightarrow & -\frac{1}{y^1_{xx}}d{\xi}^2 - \frac{1}{y^1_{xx}} {\xi}^1 = {\eta}^1  \\ 
{\mu}^2 & \rightarrow & \frac{1}{y^2_{xx}}{\xi}^1 = {\eta}^2  \\
{\mu}^3 & \rightarrow & - \frac{y^1_x}{y^1_{xx}}d{\xi}^2 + (\frac{y^1_x}{y^1_{xx}} - \frac{y^2_x}{y^2_{xx}})
{\xi}^1 + {\xi}^2 = {\eta}^3
\end{array}  \right.  \] 
This parametrization is injective because ${\xi}^1=y^2_{xx} {\eta}^2, \,\,\,\,       {\xi}^2= {\eta}^3 + 
y^2_x {\eta}^2  - y^1_x {\eta}^1$. Hence, we can replace {\it any} solution $({\eta}^1,{\eta}^2,{\eta}^3)$ of the linearized system by {\it any} couple $({\xi}^1,{\xi}^2)$, a result not evident at first sight.  \\
Finally, considering the two parametrization vertical $1$-forms:  \\
\[  {\omega}^1= \frac{1}{y^2_{xx}}{\xi}^1= {\eta}^2=\delta y^2, \hspace{1cm}  
       {\omega}^2 = \delta y^3 + y^2_x \delta y^2 - y^1_x \delta y^1   \]
we have:  \\
\[  \delta {\omega}^1=0 \Rightarrow {\omega}^1 \wedge {\omega}^2 \wedge \delta {\omega}^1=0, \hspace{1cm}   {\omega}^1 \wedge {\omega}^2 \wedge \delta {\omega}^2= \delta y^1 \wedge \delta y^2 \wedge \delta y^3 \wedge  \delta y^1_x \neq 0   \] 
and cannot therefore use the Frobenius theorem in order to integrate this vertical exterior system.    \\
According to what has been said, the linear and the nonlinear systems are both controllable. In particular, if the nonlinear system should {\it not} be controllable, it means that there should exists at least one autonomous element in $L$ that should be constrained by at least one OD equation. The linearization of such an element should produce a torsion element in $L$. The striking feature of this example is that  one can prove that $L$ is a purely differentially transcendental extension of the ground field $K$. Indeed, we may rewrite the system like:  \\
\[       d(2y^3 - (y^1 - y^2)d(y^1 + y^2)) + (y^1 - y^2) d^2 (y^1 + y^2) = 0  \]
Setting:  \\
\[          z^1 =  2 y^3 -(y^1 - y^2)(y^1_x + y^2_x),  \hspace{1cm}   z^2= y^1 + y^2   \hspace{1cm}
      \Rightarrow     \hspace{1cm}              - \frac{z^1_x}{z^2_{xx}}  = y^1 -  y^2 \]
we obtain the second order {\it nonlinear} parametrization:   \\
\[    2 \,y^1= - \frac{z^1_x}{z^2_{xx}} + z^2, \hspace{15mm} 2\,y^2= \frac{z^1_x}{z^2_{xx}}+ z^2, 
           \hspace{15mm}   2\,y^3= - \frac{z^1_x z^2_x}{z^2_{xx}} + z^1      \]
and thus $L=K<z^1,z^2>$. Hence, introducing ${\zeta}^1=\delta z^1, {\zeta}^2=\delta z^2$, we can similarly replace {\it any} solution $({\eta}^1,{\eta}^2,{\eta}^3)$ of the linearized system by {\it any} couple $({\zeta}^1,{\zeta}^2)$, a result even less evident at first sight . The new parametrization is also injective but, contrary to the previous situation, we have now $ {\bar{\omega}}= \zeta = \delta z \Rightarrow \delta {\bar{\omega}}=0 $. As a (difficult) exercise of formal integrability, we let the reader prove that the second order $2 \times 2$ operator matrix $\xi \rightarrow \zeta$ is an isomorphism (See [21], p 821 for more details on this substitution).   \\

\noindent
{\bf EXAMPLE 3.3.10}: With $n=2, m= 3, q=1, K= \mathbb{Q}$, let us consider the differential ideal $\mathfrak{p}$ generated by the two differential polynomials:    \\
\[   \hspace{2cm} P_1 \equiv y^1_2 -y^3 y^1_1 =0, \hspace{2cm}   P_2 \equiv y^2_2 - y^3 y^2_1=0    \] 
The corresponding system is easily sen to be involutive and we have:   \\
\[   K\{y^1,y^2,y^3\}/\mathfrak{p}= K\{y^3\}[y^1,y^1_1, ...;y^2, y^2_1,...]  \] 
as an integral domain and $\mathfrak{p}$ is thus a prime differential ideal allowing to define the differential field $L=Q(K\{y\}/\mathfrak{p})$ with differential transcendence basis $K<y^3> \subset L$. \\
The linearized system ${\cal{D}}_1\eta=0$ is:  \\
\[    d _2{\eta}^1 - y^3 d_1{\eta}^1 - y^1_1 {\eta}^3=0,  \hspace{2cm}
      d_2{\eta}^2 - y^3 d_1{\eta}^2 - y^2_1 {\eta}^3 = 0   \]
Multiplying the first by $ {\lambda}^1$, the second by ${\lambda}^2$ and integrating by parts, we obtain the adjoint operator $ad({\cal{D}}_1)\lambda = \mu$:  \\               
\[  \left \{  \begin{array}{rcl}                       
{\eta}^1& \rightarrow  &  - d_2 {\lambda}^1 + y^3 d_1 {\lambda}^1+ y^3_1 {\lambda}  ^1= {\mu}^1  \\
{\eta}^2& \rightarrow & - d_2  {\lambda}^2+y^3 d_1 {\lambda}^2 + y^3_1 {\lambda}^2 = {\mu}^2   \\
 {\eta}^3 &  \rightarrow &    - y^1_1 {\lambda} ^1 - y^2_1{\lambda}^2= {\mu}^3                                                                                   
 \end{array}   \right.    \]
However, this operator is {\it not} involutive because it is not even formally integrable. Nevertheless, adding the first order PD equation obtained by prolonging the zero order equation with  respect to $x^1$, we obtain:  \\
\[   \left \{  \begin{array}{rcl}
{\eta}^1& \rightarrow  &  - d_2 {\lambda}^1 + y^3 d_1 {\lambda}^1+ y^3_1 {\lambda}  ^1= {\mu}^1  \\
{\eta}^2& \rightarrow & - d_2  {\lambda}^2+y^3 d_1 {\lambda}^2 + y^3_1 {\lambda}^2 = {\mu}^2   \\
&   &   - y^1_1 d_1{\lambda}^1 - y^2_1d_1 {\lambda}^2 - y^1_{11} {\lambda}^1 - y^2_{11}{\lambda}^2=
d_1{\mu}^3   \\
 {\eta}^3 &  \rightarrow &    - y^1_1 {\lambda}^1 - y^2_1{\lambda}^2= {\mu}^3                                                                                   
 \end{array}   \right.  \fbox { $  \begin{array}{cc}
 1 & 2 \\
 1 & 2  \\
1 & \bullet  \\
\bullet & \bullet
\end{array}  $ }  \]
We obtain the unique generating first order CC $ad({\cal{D}})\mu=0$, namely:      \\
\[   d_2{\mu}^3 - y^3 d_1 {\mu}^3 - y^1_1 {\mu}^1 - y^2_1 {\mu}^2 - 2 y^3_1 {\mu}^3 =0  \]
Multiplying this CC by the test function $\xi$ and integrating by parts, we get ${\cal{D}}\xi=\eta$ over 
$L$: \\
\[    \left \{  \begin{array} {rcl}
{\mu}^1 & \rightarrow & -  y^1_1 \xi = {\eta}^1  \\
{\mu}^2 &  \rightarrow & -  y^2_1 \xi = {\eta}^2  \\
{\mu}^3 & \rightarrow & -  d_2 \xi + y^3 d_1 \xi - y^3_1\xi = {\eta}^3                          
\end{array}  \right.   \]
Sustituting, we check the two CC described by ${\cal{D}}_1 \eta=0$ plus an additional zero order CC providing the torsion element $\omega= y^1_1 \delta y^2 -y^2_1 \delta y^1$ generating $ext^1(M)$ as we have indeed $d_2 \omega - y^3 d_1 \omega - y^3_1 \omega= 0$. We let the reader check that  $\omega \wedge \delta \omega \neq 0$. It follows that $\omega$ cannot have any integrating factor according to the Frobenius Theorem.    \\ 
Eliminating $y^3$, we are left with $(y^1,y^2)$ and the nonlinear system 
$y^1_1 y^2_2 - y^1_2 y^2_1=0$ with the same conlusions as before. On the contrary, we let the reader check that the module of K\"{a}hler differentials is torsion-free for the nonlinear system 
$y^1_1 y^2_2 - y^1_2 y^2_1=1$.  \\   \\

\noindent
{\bf 4) CONCLUSION}  \\

The author of this paper got his PhD thesis under the supervising of Prof. A. Lichnerowicz and has been collaborating with him till his death in 1998. At the end of his life, he became more and more convinced that the variational origin of mathematical physics (elasticity, electromagnetism, general relativity) through the corresponding Euler-Lagrange equations was a kind of "{\it screen} " hiding a more important concept allowing to describe the "duality" existing between "fields " and "inductions ". After the author  discovered in 1995 the impossibility to parametrize the Einstein operator along a challenge proposed by J. Wheeler in 1970, he did notice that, in control theory,  "{\it a control system is controllable if and only if it is parametrizable} " and that the "screen " is just the "{\it differential double duality} " involved in differential homological algebra through the use of the "{\it extension modules} " (See Zbl 1079.93001). Hence it remained to study the systematic use of the {\it formal adjoint} in the noncommutative situation met when linearizing nonlinear systems of OD or PD equations. Among the best useful examples, one has the following {\it differential sequence}, indicating below the fiber dimensions of the vector bundles involved with $F_0=S_2T^*$:  \\
\[  0 \rightarrow \Theta \rightarrow T \stackrel{Killing}{\longrightarrow} F_0 \stackrel{Riemann}{\longrightarrow} F_1 \stackrel{Bianchi}{\longrightarrow} F_2 
\rightarrow ...  \]
\[0 \rightarrow \Theta \rightarrow n \rightarrow n(n+1)/2 \rightarrow n^2(n^2-1)/12 \rightarrow n^2(n^2-1)(n-2)/24 \rightarrow ...  \]
where $\Theta$ is the sheaf of Killing vector fields for the Euclidean metric when $n=3$ in elasticity or the Minkowski metric when $n=4$ in general relativity. Defining the adjoint operators:  \\  
\[       Cauchy=ad(Killing), \,\, Beltrami=ad(Riemann), \,\, Lanczos=ad(Bianchi)   \] 
one discovers that Lanczos was in fact dreaming to construct the {\it adjoint differential sequence}: \\
\[  0 \leftarrow ad(T) \stackrel{Cauchy}{\longleftarrow} ad(F_0) \stackrel {Beltrami}{\longleftarrow} ad(F_1)\stackrel{Lanczos}{\longleftarrow} ad(F_2) \leftarrow ...   \]
where $ad(E)={\wedge}^nT^*\otimes E^*$ for any vector bundle $E$ where $E^*$ is obtained from $E$ by inverting the transition rules when changing local coordinates. Accordingly, the only true problem left was to prove that {\it each operator is indeed parametrized by the preceding one} in both sequences, a highly non evident fact. Similarly, we have 
the Poincar\'{e} sequence for the exterior derivative:  \\
\[       T^* \stackrel{d}{\longrightarrow} {\wedge}^2 T^* \stackrel{d}{\longrightarrow}{\wedge}^3T^*  \]
\[     \hspace{15mm}  n \rightarrow n(n-1)/2 \rightarrow n(n-1)(n-2)/6    \]
with $dA=F \Rightarrow dF=0$ for electromagnetism when $n=4$ where $A$ is the EM potential and $F$ the EM field describing the first Maxwell operator and its parametrization. The adjoint sequence:  \\
\[       ad(T^*) \stackrel{ad(d)}{\longleftarrow}  ad({\wedge}^2T^*) \stackrel{ad(d)}{\longleftarrow} ad({\wedge}^3T^*)   \]
is used for the EM induction and second Maxwell operator both with its parametrization by the so-called pseudo-potential. In both situations, there is no need to appeal to variational calculus which is only used for exhibiting the respective constitutive laws. We hope that the many tricky examples presented in this paper will be used later on as test-examples for computer algebra.  \\    \\

\noindent
{\bf REFERENCES   }  \\ 
    
\noindent
[1] Aksteiner, S., Andersson, L., B\"{a}ckdahl, T.: New Identities for Linearized Gravity on the Kerr Space-time, https://arxiv.org/1601.06084v3  \\
\noindent
[2] Bialynicki-Birula, A.: On Galois Theory of Fields with Operators, Amer.J. Math., 84 (1962) 89-109.  \\
\noindent
[3] Bjork, J.E.: Analytic D-Modules and Applications, Kluwer (1993).  \\ 
\noindent
[4] Bourbaki, N.: Alg\`{e}bre, Chapter 10, Alg\`{e}bre Homologique, Springer, 2006.  \\
\noindent
[5] Cosserat, E., \& Cosserat, F.: Th\'{e}orie des Corps D\'{e}formables, Hermann, Paris, (1909).\\
\noindent
[6] Janet, M.: Sur les Syst\`{e}mes aux D\'{e}riv\'{e}es Partielles, Journal de Math., 8 (1920) 65-151. \\
\noindent 
 \noindent
[7] Johnson, J.: K\"{a}hler Differentials and Differential Algebra, Trans. Amer. Math. Society, 192 (1974) 201- 208.  \\
\noindent
[8] Johnson, J.: Prolongations of Integral Domains, Journal of Algebra, 94 (1985) 173-210.  \\ 
\noindent
[9] Kaplansky, I.: An Introduction to Differential Algebra, Hermann, Paris ((1957, 1976).   \\
\noindent
[10] Kashiwara, M.: Algebraic Study of Systems of Partial Differential Equations, M\'{e}moires de la Soci\'{e}t\'{e} Math\'{e}matique de France, 63 (1995) (Transl. from Japanese of his 1970 MasterÕs Thesis).  \\
\noindent
[11] Kolchin, E. R.: Differential Algebra and Algebraic Groups, Academic Presss, (1973).  \\
\noindent
[12] Macaulay, F.S.: The Algebraic Theory of Modular Systems, Cambridge Tract 19, Cambridge University Press, London, 1916 (Reprinted by Stechert-Hafner Service Agency, New York, 1964).  \\
\noindent
[13] Northcott, D.G.: An Introduction to Homological Algebra, Cambridge university Press (1966).  \\
\noindent
[14] Oberst, U.: Multidimensional Constant Linear Systems, Acta Appl. Math., 20 (1990) 1-175.  \\
\noindent
[15] Palamodov, V.P.: Linear Differential Operators with Constant Coefficients, Grundlehren der Mathematischen Wissenschaften, 168, Springer-Verlag (1970).  \\ 
\noindent
[16] Pirani, F.A.E., Robinson, D.C., Shadwick, W.F.: Local Jet Bundle Formulation of B\"{a}cklund Transformations, MPST 1, D. Deidel, Springer (1979).  \\
\noindent
[17] Pommaret, J.-F.: Systems of Partial Differential Equations and Lie Pseudogroups, Gordon and Breach, New York (1978); Russian translation: MIR, Moscow,(1983).\\
\noindent
[18] Pommaret, J.-F.: Differential Galois Theory, Gordon and Breach, New York (1983).\\
\noindent
[19] Pommaret, J.-F.: Lie Pseudogroups and Mechanics, Gordon and Breach, New York (1988).\\
\noindent
[20] Pommaret, J.-F.: Partial Differential Equations and Group Theory, Kluwer (1994).\\
http://dx.doi.org/10.1007/978-94-017-2539-2    \\
\noindent
[21] Pommaret, J.-F.: Partial Differential Control Theory, Kluwer, Dordrecht (2001) (Zbl 1079.93001).   \\
\noindent
[22] Pommaret, J.-F.: Algebraic Analysis of Control Systems Defined by Partial Differential Equations, in "Advanced Topics in Control Systems Theory", Springer, Lecture Notes in Control and Information Sciences 311 (2005) Chapter 5, pp. 155-223.\\
\noindent
[23] Pommaret, J.-F.: Parametrization of Cosserat Equations, Acta Mechanica, 215 (2010) 43-55.\\
http://dx.doi.org/10.1007/s00707-010-0292-y  \\
\noindent
[24] Pommaret, J.-F.: Spencer Operator and Applications: From Continuum Mechanics to Mathematical Physics, in "Continuum Mechanics-Progress in Fundamentals and Engineering Applications", Dr. Yong Gan (Ed.), ISBN: 978-953-51-0447--6, InTech (2012) Available from: \\
http://dx.doi.org/10.5772/35607   \\
\noindent
[25] Pommaret, J.-F.: The Mathematical Foundations of General Relativity Revisited, Journal of Modern Physics, 4 (2013) 223-239. \\
 https://dx.doi.org/10.4236/jmp.2013.48A022   \\
  \noindent
[26] Pommaret, J.-F.: Macaulay Inverse Systems and Cartan-K\"ahler Theorem,  \\
 https://arxiv.org/abs/1411.7070   \\
\noindent
[27] Pommaret, J.-F.: Relative Parametrization of Linear Multidimensional Systems, Multidim. Syst. Sign. Process., 26 (2015) 405-437.  \\
DOI 10.1007/s11045-013-0265-0   \\
\noindent
[28] Pommaret, J.-F.: Airy, Beltrami, Maxwell, Einstein and Lanczos Potentials Revisited, Journal of Modern Physics, 7 (2016) 699-728. \\
\noindent
https://dx.doi.org/10.4236/jmp.2016.77068   \\
\noindent
[29] Pommaret, J.-F.: Deformation Theory of Algebraic and Geometric Structures, Lambert Academic Publisher (LAP), Saarbrucken, Germany (2016). A short summary can be found in "Topics in Invariant Theory ", S\'{e}minaire P. Dubreil/M.-P. Malliavin, Springer 
Lecture Notes in Mathematics, 1478 (1990) 244-254.\\
https://arxiv.org/abs/1207.1964  \\
\noindent
[30] Pommaret, J.-F.: New Mathematical Methods for Physics, Mathematical Physics Books, Nova Science Publishers, New York (2018) 150 pp.  \\
\noindent
[31]  Pommaret, J.-F.: The Mathematical Foundations of Elasticity and Electromagnetism Revisited, Journal of Modern Physics, 10 (2019) 1566-1595.     \\
 https://doi.org/10.4236/jmp.2019.1013104 (https://arxiv.org/abs/1802.02430 ) \\
\noindent
[32] Pommaret, J.-F.: Generating Compatibility Conditions and General Relativity, J. of Modern Physics, 10, 3 (2019) 371-401.  \\
\noindent
https://doi.org/10.4236/jmp.2019.103025   \\
\noindent
[33] Pommaret, J.-F.: Differential Homological Algebra and General Relativity, J. of Modern Physics, 10 (2019) 1454-1486. \\
\noindent
https://doi.org/10.4236/jmp.2019.1012097   \\
\noindent
[34] Pommaret, J.-F.: The Conformal Group Revisited, https://arxiv.org/abs/2006.03449   \\
\noindent
[35]  Pommaret, J.-F.: A Mathematical Comparison of the Schwarzschild and Kerr Metrics, Journal of modern Physics, 11, (2020) 1672-1710.  
https://dx.doi.org/10.4236/jmp.2020.1110104     \\
\noindent
[36] Pommaret, J.-F.: Minimum Parametrization of the Cauchy Stress Operator, Journal of Modern Physics, 12 (2021) 453-482. https://doi.org/10.4236/jmp.2021.124032  \\
\noindent
[37] Pommaret, J.-F.: Homological Solution of the Lanczos Problems in Arbitrary Dimension, Journal of Modern Physics, 12 (2021) 829-858. 
https://doi.org/10.4236/jmp.2021.126053  \\
\noindent
[38]Quadrat, A.: Une Introduction \`{a} l'Analyse Alg\'{e}brique Constructive et \`{a} ses Applications, INRIA Research Report 7354, AT-SOP Project, july 2010. 
Les Cours du CIRM, 1 no. 2: Journ\'{e}es Nationales de Calcul Formel (2010), p281-471 (doi:10.5802/ccirm.11). \\
\noindent 
[39] Rotman, J.J.: An Introduction to Homological Algebra, Pure and Applied Mathematics, Academic Press (1979).  \\
\noindent
[40] Schneiders, J.-P.: An Introduction to D-Modules, Bull. Soc. Roy. Sci. Li\`{e}ge, 63, 223-295 (1994).  \\
\noindent
[41] Shankar, S.: Controllability and Vector Potential, https://arxiv.org/abs/1911.01238v2  \\
\noindent
[42] Shankar, S.: The Canonical Controller for Distributed Systems, Multidimensional Systems and Signals Processing, 32 (2021) 303-311. 
https://arxiv.org/10.1007/s11045-020-00740-1    \\
\noindent
[43] Spencer, D.C.: Overdetermined Systems of Partial Differential Equations, Bull. Am. Math. Soc., 75 (1965) 1-114.\\
\noindent
[44] Zerz, E.: Topics in Multidimensional Linear Systems Theory, Lecture Notes in Control and Information Sciences, Springer, LNCIS 256 (2000). \\
\noindent
[45] Zerz, E.: An Algebraic Analysis Approach to Linear Time-Varying Systems, IMA Journal of Mathematical Control and Information, 23 (2005) 113-126. \\
https://doi.org/10.3182/20050703-6-CZ-1902.00200  \\

\end{document}